\documentclass[11pt]{amsart}
\usepackage{amsmath}
\usepackage{amssymb}
\usepackage[a4paper,left=3cm,right=3cm]{geometry}
\usepackage{xcolor}
\usepackage{multirow}
\usepackage{stackengine}
\usepackage{scalerel}
\usepackage{mathtools}
\usepackage{enumitem}
\usepackage[unicode=true,pdfusetitle, bookmarks=true,bookmarksnumbered=false,bookmarksopen=false,breaklinks=false,pdfborder={0 0 0},pdfborderstyle={},colorlinks=true]{hyperref}

\hypersetup{linkcolor={blue!50!black},citecolor={green!60!black},urlcolor={blue!80!black}}

\setlist[enumerate,1]{label={(\arabic*)}}

\newtheorem{theorem}{Theorem}[section] %
\newtheorem{corollary}[theorem]{Corollary}
\newtheorem{lemma}[theorem]{Lemma}
\newtheorem{proposition}[theorem]{Proposition}

\theoremstyle{definition}
\newtheorem{definition}[theorem]{Definition}
\newtheorem{example}[theorem]{Example}
\newtheorem{remark}[theorem]{Remark}
\newtheorem{fact}[theorem]{Fact}
\newtheorem{notation}[theorem]{Notation}

\newcommand{\Tref}[1]{{Theorem~\ref{#1}}}
\newcommand{\Lref}[1]{{Lemma~\ref{#1}}}
\newcommand{\Pref}[1]{{Proposition~\ref{#1}}}
\newcommand{\Cref}[1]{{Corollary~\ref{#1}}}
\newcommand{\Dref}[1]{{Definition~\ref{#1}}}
\newcommand{\Exref}[1]{{Example~\ref{#1}}}
\newcommand{\Sref}[1]{{Section~\ref{#1}}}

\newcommand\suchthat{\,\ifnum\currentgrouptype=16\middle\fi|\,}

\newcommand{\mdash}
           {\discretionary{}{}{\kern 0.16667em}---\discretionary{}{}{\kern 0.16667em}}
\DeclareMathOperator{\rank}{rk}
\DeclareMathOperator{\Rank}{rank}
\DeclareMathOperator{\ann}{ann}
\DeclareMathOperator{\Span}{Span}
\DeclareMathOperator{\Image}{Im}
\DeclareMathOperator{\Lie}{Lie}
\DeclareMathOperator{\GK}{GKdim}
\DeclareMathOperator{\Free}{Fr}
\DeclareMathOperator{\diag}{diag}
\newcommand{\set}[1]{{\left\{#1\right\}}}
\newcommand\M[1][n]{{\operatorname{M}_{#1}}}
\newcommand\GL[1][n]{{\operatorname{GL}_{#1}}}
\newcommand{\eps}{\varepsilon}
\newcommand{\Id}{\mathrm{Id}} 
\newcommand{\sub}{\subseteq}
\newcommand{\C}{\mathbb{C}}
\newcommand{\N}{\mathbb{N}}
\newcommand{\Alg}{\mathcal{A}}
\newcommand{\dist}{\operatorname{dist}}

\begin{document}

\title{Rank-Stability of Polynomial Equations}

\author{Tomer Bauer}
\address{Department of Mathematics, Bar-Ilan University, Ramat Gan 5290002, Israel.}
\email{mathzeta2@gmail.com}

\author{Guy Blachar}
\address{Department of Mathematics, Weizmann Institute of Science, Rehovot 7610001, Israel.}
\email{guy.blachar@gmail.com}

\author{Be'eri Greenfeld}
\address{Department of Mathematics, University of Washington, Seattle, WA 98195-4350, USA.}
\email{grnfld@uw.edu}

\begin{abstract}
Extending the thoroughly studied theory of group stability, we study Ulam stability type problems for associative and Lie algebras. Namely, we investigate obstacles to rank-approximation of ``almost'' solutions by exact solutions to systems of polynomial equations. This leads to a rich theory of stable associative and Lie algebras, with connections to linear soficity, amenability, growth, and group stability. We develop rank-stability and instability criteria, examine the effect of algebraic constructions on rank-stability, and prove that while finite-dimensional associative algebras are rank-stable, ``most'' finite-dimensional Lie algebras are not.
\end{abstract}

\maketitle

\section{Introduction}

Ulam's stability problem asks when an approximate solution to a given equation, with respect to some metric, is ``close'' to an exact solution, see \cite{Ulam60}. 
To make this well-defined, we need to specify the equations and the distance function. 
For instance, a famous case of this problem is due to Halmos \cite{Halmos76}, who asked whether self-adjoint matrices that are almost commuting are close to commuting self-adjoint matrices, where ``almost commuting'' and ``close'' refer to the operator norm. This was solved in the affirmative by Lin \cite{Lin96} (see also \cite{FriisRordam96}). Interestingly, a similar result for triples of matrices is false \cite{Davidson85}.
Voiculescu \cite{Voiculescu83} showed that there exists a sequence of pairs of unitary matrices, for which the commutators approach to the identity matrix with respect to the operator norm, but for which there are no arbitrarily good approximations by commuting unitary matrices (see also \cite{ExelLoring89}). For other variations of this problem and related results, see \cite{Glebsky10,Hadwin98,Rosenthal69,VonNeumann29}. 
A discrete version of the stability problem for the commutativity equation $xy=yx$ has been resolved by Arzhantseva and P\u{a}unescu \cite{ArzhantsevaPaunescu15}, who proved that permutations that almost commute with respect to the normalized Hamming distance are close to commuting permutations.

Stability of more general equations has been systematically studied by Glebsky and Rivera \cite{GlebskyRivera09}.
Given a finite set of equations in a free group $\Free(x_1,\dots,x_d)$, say the equations $w_1(\vec{x})=1,\dots,w_r(\vec{x})=1$, we say that the quotient group $\left<x_1,\dots,x_d\ |\ w_1,\dots,w_r \right>$ is \textbf{stable} with respect to a sequence of groups with distance functions $\{(G_n,\dist_n)\}_{n=1}^{\infty}$ if for every $\varepsilon>0$ there exists $\delta>0$ such that for every tuple of elements $g_1,\dots,g_d\in G_n$ satisfying $\dist_n(w_i(\vec{g}),1)<\delta$ for all $1\leq i\leq r$, there exist $g'_1,\dots,g'_d\in G_n$ for which $w_1(\vec{g'})=\cdots=w_r(\vec{g'})=1$ and $\dist_n(g_j,g'_j)<\varepsilon$ for all $1\leq j\leq d$. 
In particular, stability of the commutation equation is just stability of abelian groups; finite groups are stable with respect to the Hamming distance on symmetric groups \cite{GlebskyRivera09} (this is called P-stability).
A deep connection of group stability with vanishing cohomology has been found in \cite{ChiffreGlebskyLubotzkyThom20}.  
A remarkable result of Becker, Lubotzky and Thom \cite{BeckerLubotzkyThom19} characterizes P-stability of amenable groups by means of invariant random subgroups, and proves that virtually-polycyclic groups are P-stable. Stability results of groups have connection to property testing \cite{BeckerLubotzkyMosheiff23} and error-correcting codes \cite{ChapmanVidickYuen23}.
For other recent interesting results, see \cite{LazarovichLevit23,LazarovichLevitMinsky19}.

In the context of Lie algebras and non-commutative ($C^{*}$-)algebras, approximate representations were studied by Ioos, Kazhdan and Polterovitch \cite{IoosKazhdanPolterovich23} with interesting connections to geometric quantizations of spheres and tori. In a similar spirit, soficity of associative algebras was studied by Arzhantseva and P\u{a}unescu \cite{ArzhantsevaPaunescu17} (and by Cinel \cite{Cinel23} in the context of Lie algebras), where the role of the normalized Hamming distance is played by the normalized rank distance; see also \cite{Elek05}.

In this vein, Elek and Grabowski \cite{ElekGrabowski21} proved that tuples of almost-commuting unitary matrices (with respect to the normalized rank distance) are close to tuples of genuinely commuting matrices.

\subsection{Rank-stability}
In this paper we study stability of algebras with respect to the normalized rank distance. This is a natural framework that enables one to study rank-stability of non-commutative polynomial matrix equations, of Lie algebras (via their universal enveloping algebras) and of groups (via their group algebras).

For a matrix $A\in\M[n](F)$ over a field $F$, let its \textbf{rank} be $\Rank(A)=\dim_F\Image A$ and
its \textbf{normalized rank} be $\rank(A)=\frac{1}{n} \Rank(A)$.
This endows $\M[n](F)$ with a discrete metric by $d(A,B)=\rank(A-B)$.
Following \cite{ElekGrabowski21}, let $\widehat{A}\in \M[\N](F)=\bigcup_{n=1}^{\infty} \M[n](F)$ be the finitely supported matrix obtained from $A$ by setting $\widehat{A}_{i,j}=A_{i,j}$ if $1\le i,j\le n$, and $\widehat{A}_{i,j}=0$ otherwise. Notice that $\Rank(\widehat{A})=\dim_F \Image \widehat{A} = \Rank(A)$ is well-defined, and that
$\widehat{\cdot}$ is a (non-unital) embedding $\M[n](F)\hookrightarrow \M[\N](F)$.

Given a system of (non-commutative) polynomials $P_1,\dots,P_r$ in $d$ variables, we say that a tuple $\vec{B}=(B_1,\dots,B_d)$ of matrices in $\M[n](F)$ is a \textbf{solution} of the system if
\[
    P_i(B_1,\dots,B_d) = 0
\]
for all $1\le i\le r$. Scalar monomials participating in $P_i$ are interpreted as scalar $n\times n$ matrices in $P_i(\vec{B})$. For emphasis, we sometimes say that $\vec{B}$ is an \textbf{exact solution} of the system, in contrast to an approximate solution (or approximate representation) that is at the heart of the following definition, which is the algebraic counterpart of group rank-stability defined in \cite{ElekGrabowski21}.

\begin{definition}[{Rank-stability of algebras}]\label{def:rank-stab}
Let $\Alg=F\left<x_1,\dots,x_d\right>/\left<P_1,\dots,P_r\right>$ be a finitely presented algebra. We say that $\Alg$ is \textbf{rank-stable} if for every $\eps>0$ there exists $\delta>0$ such that for every $n\in\N$ and every $d$-tuple of $n\times n$ matrices $\vec{A}=(A_1,\dots,A_d)$ satisfying:
\[ \rank P_i(A_1,\dots,A_d) <
\delta \ \ \text{for all}\ 1\le i\le r
\]
there exists $n'\in\N$ and a solution $\vec{B}=(B_1,\dots,B_d)$ of $n'\times n'$ matrices of $P_1,\dots,P_r$ such that
\[
\Rank(\widehat{A_i}-\widehat{B_i}) < \eps n \ \ \text{for all}\ 1\le i\le d.
\]
We say that $\vec{B}=(B_1,\dots,B_d)$ \textbf{$\eps$-approximates} $\vec{A}=(A_1,\dots,A_d)$.
\end{definition}

The use of $\,\widehat{\cdot}\,$ and the consideration of different matrix dimensions in the approximating solutions, is necessary in order to avoid trivial obstructions, and will be made clear in Example \ref{ex:mat}. 
\Dref{def:rank-stab} is given in \cite{ElekGrabowski21} in the group-theoretic context, where the authors prove that $\mathbb{Z}^k$ is rank-stable when all matrices considered are unitary or self-adjoint. The question of whether $\mathbb{Z}^k$ is rank-stable without any spectral restrictions on the involved matrices remains open \cite[Remark~2]{ElekGrabowski21}.

\Dref{def:rank-stab} is formulated for associative algebras but applies mutatis-mutandis to Lie algebras, in which case $F\left<x_1,\dots,x_d\right>$ is replaced by the finitely generated free Lie algebra $\Lie\left<x_1,\dots,x_d\right>$, with $\M[n](F)$ thought of as $\mathfrak{gl}_n(F)$ and $P_1,\dots,P_r$ are Lie polynomials. For instance, every finite-dimensional Lie algebra with basis $\{e_1,\dots,e_d\}$ can be presented as:
\[
\Lie\left<e_1,\dots,e_d\right>/\left<[e_i,e_j]-{\textstyle \sum_{k=1}^d} \gamma_{ij}^k e_k\ \suchthat\ 1\leq i<j\leq d \right>
\]
for some $\gamma_{ij}^k \in F$. 
A finitely presented Lie algebra $\mathfrak{g}$ 
 is rank-stable if and only if its universal enveloping algebra $U(\mathfrak{g})$ is rank-stable. Indeed, given a Lie algebra $\mathfrak{g}$ with presentation $\Lie\left<x_1,\dots,x_d\right>/\left<P_1,\dots,P_r\right>$, its universal enveloping algebra is presented by $F\left<x_1,\dots,x_d\right>/\left<P_1,\dots,P_r\right>$ where each $P_i$ is interpreted as an associative polynomial, namely $[a,b]=ab-ba$, and the notions of $\rank$ and $\varepsilon$-approximation are identical in the associative and Lie settings.

\subsection{Results and structure of the paper}

In \Sref{sec:basic} we give some basic properties and examples of rank-stable algebras, which support the main definition as well as the role of the existence of finite-dimensional representations; we also discuss the connection with linear soficity.
In \Sref{sec:compression} we develop a useful ``compression machinery'' which allows us to effectively control the dimensions of exact solutions of equations by means of approximate solutions.
In \Sref{sec:independence} we prove that rank-stability is indeed an algebra property rather than a property of a given presentation (\Tref{thm:two-pres}).

In \Sref{sec:finite-dimensional} we prove that finite-dimensional associative algebras are rank-stable (\Tref{thm:fin-dim}); this is parallel to the group-theoretic counterpart, but differs from the Lie-theoretic situation.
\Sref{sec:group} proves that rank-stability of a group is equivalent to the rank-stability of its group algebra (\Tref{thm:group-alg}).
\Sref{sec:constructions} shows that rank-stability is preserved under free products, direct products and matrix rings.

Finally, in \Sref{sec:instability} we produce a non-stability machinery that demonstrates the ubiquity of non rank-stable algebras (\Tref{thm:unstable}), even within classes of tamely behaved algebras, and derive applications to Lie algebras (\Tref{thm:Lie}).

\smallskip
\noindent \textit{Notation and terminology.} 
Given matrices $A_1,\dots,A_d\in\M[n](F)$, we denote for convenience $\vec{A}=(A_1,\dots,A_d)$.
For a polynomial $f(x_1,\dots,x_d)$, we denote by $\widehat{f}(\vec{A})=\widehat{f}(A_1,\dots,A_d)\in\M[\N](F)$ the matrix obtained by applying $\,\widehat{\cdot}\,$ to $f(\vec{A})$.
All polynomials in this paper are assumed to be polynomials over non-commutative variables, unless explicitly stated.

\section{Basic Properties and Examples} \label{sec:basic}

We will liberally use the following elementary proposition in proofs and calculations without mentioning it, and introduce some common notation.

\begin{proposition}[{cf.~\cite[Proposition~2.2]{ArzhantsevaPaunescu17}}]
Let $A, A'\in\M[n](F)$ and $B\in\M[m](F)$. Then:
\begin{enumerate}
    \item $\Rank A = 0$ if and only if $A=0_n$;
    \item $\Rank(\Id_n) = n$, where $\Id_n\in\M[n](F)$ is the identity matrix;
    \item $\Rank(A+A') \le \Rank A + \Rank A'$;
    \item $\Rank(AA') \le \min\set{\Rank A, \Rank A'}$;
    \item $\Rank(A^{-1}A'A) = \Rank A'$ if $A$ is invertible;
    \item $\Rank(A \oplus B)=\Rank A +\Rank B$, where $A \oplus B$ is a block diagonal matrix;
    \item $\Rank(A \otimes B)=\Rank A \cdot\Rank B$, where $A \otimes B$ is the Kronecker product;
    \item Suppose that $B$ is the matrix obtained from $A$ by adding zero rows and columns if $n\le m$, or by removing the last $n-m$ rows and columns if $n>m$. Then $\Rank(\widehat{A}-\widehat{B})\le 2|n-m|$.

\end{enumerate}
\end{proposition}

The notion of stability is closely related to soficity (cf.~\cite{BeckerLubotzkyThom19,ElekGrabowski21}). Recall the notion of linear soficity from~\cite{ArzhantsevaPaunescu17}:
let $\{n_k\}_{k=1}^{\infty}$ be an infinite set of positive integers and let $\mathcal{U}$ be a non-principal ultrafilter on $\N$. Let $\rho_{\mathcal{U}}\colon \prod_k \M[n_k](F) \to [0,1]$ be the asymptotic rank function, defined by
\[
\rho_{\mathcal{U}}(A_1,A_2,\dots) = \lim_{k\rightarrow \mathcal{U}} \rank(A_k).
\]
A countably generated algebra is \textbf{linear sofic} if it embeds into a metric ultraproduct $\prod_{k\rightarrow \mathcal{U}} \M[n_k](F) / \ker(\rho_\mathcal{U})$.
The following is in the spirit of~\cite[Remark~5(1)]{ElekGrabowski21}. We use results from \Sref{sec:compression}, which does not rely on the current section.

\begin{proposition}\label{prop:sofic}
Let $\Alg$ be a finitely presented, rank-stable algebra. Suppose that~$\Alg$ is linear sofic. Then $\Alg$ is residually finite-dimensional.
\end{proposition}

\begin{proof}
  Since $\Alg$ is linear sofic, we have an embedding:
  \[
    \Phi\colon \Alg \hookrightarrow \prod_{k\rightarrow \mathcal{U}} \M[n_k](F) / \ker(\rho_\mathcal{U}).
  \]
  Fix a presentation
  \[
    \Alg=F\left<x_1,\dots,x_d\right>/\left<P_1,\dots,P_r\right>
  \]
  that is rank-stable. For any $k\ge 1$, let $X_{1,k},\dots,X_{d,k}\in\M[n_k](F)$ denote the images of $x_1,\dots,x_d$ (respectively) in $\M[n_k](F)$ under $\Phi$.

  Let $0\ne a\in \Alg$, and write $a=p(x_1,\dots,x_d)$ for some $d$-variate polynomial $p$. Let $L$ and $M$ be the number of monomials and the maximal length of monomials in $p$, respectively. Since $\Phi$ is an embedding, there is some $\eps>0$ such that:
  \[
    S_0=\set{k\suchthat \rank p(X_{1,k},\dots,X_{d,k})\ge\eps}\in\mathcal{U}.
  \]
Let $\delta>0$ be the $\delta$ from the rank-stability of $\Alg$ with respect to $\eps'=\frac{\eps}{2(LM+d)}$; we assume that $\delta<(\eps' d)^2$. 
By the definition of the metric ultraproduct, we have
  \[
    S_1=\set{k\suchthat \forall j:\rank(P_j(X_{1,k},\dots,X_{d,k}))<\delta}\in\mathcal{U}.
  \]
  By the rank-stability of $\Alg$, for every $k\in S_1$ there exist $n_k'\ge 1$ and matrices $X_{1,k}',\dots,X_{d,k}'\in\M[n_k'](F)$ such that $P_j(X_{1,k}',\dots,X_{d,k}')=0$ for all $1\le j\le r$, and $\Rank(\widehat{X_{i,k}}-\widehat{X_{i,k}'})\le\eps' n_k$. By \Lref{lem:small-approx} and since $\delta<(\eps' d)^2$, we may assume $\left|n_k'-n_k\right|\le\eps' dn_k$. 
    Fix $k\in S_0\cap S_1$. By \Lref{lem:computing polyrank}, we have,
  \[
    \Rank(\widehat{p}(X_{1,k},\dots,X_{d,k})-\widehat{p}(X_{1,k}',\dots,X_{d,k}')) \le LM\eps' n_k + \eps' dn_k = \frac{\eps}{2} n_k.
  \]
  As $\Rank p(X_{1,k},\dots,X_{d,k})\ge\eps n_k$, we deduce
  \[
    \Rank p(X_{1,k}',\dots,X_{d,k}') \ge \frac{\eps}{2}n_k.
  \]
  In particular, $p(X_{1,k}',\dots,X_{d,k}')\ne 0_{n_k'}$.
   Therefore the homomorphism $f\colon \Alg\rightarrow \M[n_k'](F)$ defined by $x_i\mapsto X_{i,k}'$ carries $f(a)\ne 0$. It follows that $\Alg$ is residually finite-dimensional.
\end{proof}

Therefore, special attention is paid to the rank-stability of residually finite-dimensional algebras. Recall that a Lie algebra is residually finite-dimensional if and only if its universal enveloping algebra is residually finite-dimensional \cite{Michaelis87b,Michaelis87}.

\begin{example}
The first Weyl algebra $\Alg_1=\C\left<x,y\right>/\left<xy-yx-1\right>$ is infinite-dimensional and simple, so it admits no 
finite-dimensional representations. However, $\Alg_1$ admits approximate representations: for each $n>1$, the matrix
\[ \diag(1,\dots,1,-(n-1)) \in\M[n](\C)\]
has zero trace and therefore can be written as a commutator $[X_n,Y_n]$. Thus its normalized rank is $\rank \left( [X_n,Y_n] - \Id_n\right)=\frac{1}{n}$. 
Therefore, $\Alg_1$ is not rank-stable. An alternative argument for the non rank-stability of $\Alg_1$ follows from \Pref{prop:sofic}, since $\Alg_1$ is a linear sofic algebra (see~\cite{ArzhantsevaPaunescu17}) but not residually finite-dimensional.
\end{example}

An extreme case where rank-stability holds vacuously is when an algebra does not admit even approximate finite-dimensional representations. In other words, when there are no homomorphisms to any metric ultraproduct.

\begin{example}[{A vacuously rank-stable algebra}]\label{ex:vac-stab}
Let $\Alg=F\left<x,y,z\right>/\left<xyz-1,xzy\right>$. This is an infinite-dimensional algebra spanned by all monomials not containing any occurrences of $xyz$ or $xzy$; this can be seen using Bergman's diamond lemma~\cite{Bergman78}, since $xyz,xzy$ have no overlaps or self-overlaps. Furthermore, the algebra $\Alg$ admits no finite-dimensional representations, since the matrix images of $x,y,z$ must be invertible (as $xyz=1$) but one of them must be singular (as $xzy=0$). However, $\Alg$ is rank-stable for a vacuous reason: it has no $\delta$-approximate representations for any $\delta<\frac{1}{4}$. Indeed, if $X,Y,Z\in \M[n](F)$ for some $n$ are such that $\rank(XYZ-\Id_n)<\frac{1}{4}$ then $\rank(XYZ)>\frac{3}{4}$, so $\rank(X),\rank(Y),\rank(Z)>\frac{3}{4}$. Hence $\rank(XZY)>\frac{1}{4}$.
\end{example}

Notice that in order to avoid unpleasant obstructions in \Dref{def:rank-stab} it is necessary to allow the approximating solutions to have a different matrix size, namely, to use $\widehat{\cdot}$. Otherwise, even finite-dimensional algebras need not be rank-stable; indeed, \Tref{thm:fin-dim} shows that finite-dimensional algebras are rank-stable.

\begin{example} \label{ex:mat}
The matrix algebra $\M[k](F)$, for $k>1$, admits approximate representations that are not close (in the rank distance) to any representation within the same representation dimension. Fix a presentation by standard matrix units: \[ \M[k](F) = F\left<e_{ij}\ |\ 1\leq i,j\leq k \right> / \left<e_{ij}e_{lt}=\delta_{jl}e_{it},\;{\textstyle \sum_{i=1}^k}e_{ii}=1\ |\ 1\leq i,j,l,t\leq k \right>. \]
Set matrices \[ E_{ij}=\left(e_{ij}\otimes \Id_n\right)\oplus \Id_1\in \M[nk+1](F) \] for each $1\leq i,j\leq k$, where $\{e_{ij}\}_{1\leq i,j\leq k}$ are the standard $k\times k$ matrix units of $\M[k](F)$. Then for each defining relation $P(\vec{e})$ of $\M[k](F)$ in the above presentation we have $\rank P(\vec{E})\leq \frac{1}{nk+1}\xrightarrow{n\rightarrow \infty}0$. However, there are no proper representations $\M[k](F)\rightarrow \M[nk+1](F)$, since $k$ does not divide $nk+1$.
\end{example}

\begin{remark}
Let $\Alg$ be an algebra. When verifying its rank-stability we may always assume, given $\varepsilon>0$, that the dimension $n$ of a given approximate solution $A_1,\dots,A_d\in\M[n](F)$ is larger than any given integer $k \ge 1$.
Indeed, suppose that we choose $\delta<\frac{1}{k}$. If $k\ge n$ and
    \[
        \rank P_i(A_1,\dots,A_d)<\delta<\frac{1}{k},
    \]
    then
    \[
        \Rank(P_i(A_1,\dots,A_d))<\frac{n}{k}\le 1,
    \]
    so we must have $P_i(A_1,\dots,A_d)=0_n$. Therefore, by taking $\delta<\frac{1}{k}$, we may always assume that $n>k$.
\end{remark}

\subsection{Warm-up examples} We finish this section with a couple of quick examples. Example \ref{ex:p(x)} is a special case of \Tref{thm:fin-dim}.

\begin{example} \label{ex:p(x)}
  Let $p(x)\in F[x]$ be a non-constant polynomial that splits over $F$. Then the algebra $F[x]/{\left<\,p(x)\right>}$ is rank-stable. Although this is a special case of \Tref{thm:fin-dim}, stating that any finite-dimensional associative algebra is rank-stable, we give a direct proof.

  Write $p(x)=(x-\lambda_1)^{r_1}\cdots(x-\lambda_k)^{r_k}$ for distinct $\lambda_1,\dots,\lambda_k\in F$. Take $A\in\M[n](F)$ such that $\rank p(A)<\eps$. Take a basis $v_1,\dots,v_m$ of $U=\bigoplus_{j=1}^k\ker((A-\lambda_j\Id_n)^{r_j})$ such that each $v_i$ belongs to some $\ker((A-\lambda_j\Id_n)^{r_j})$. Note that by Sylvester's inequality
  \begin{align*}
      m = \sum_{j=1}^k\dim_F \ker((A-\lambda_j\Id_n)^{r_j}) &\ge \dim_F\ker\left(\prod_{j=1}^k(A-\lambda_j\Id_n)^{r_j}\right) \\
      &= \dim_F\ker p(A) > (1-\eps)n.
  \end{align*}
  Extend $v_1,\dots,v_m$ to a basis $v_1,\dots,v_m,v_{m+1},\dots,v_n$ of $F^n$, and define a matrix $B\in\M[n](F)$ by $Bv_i=Av_i$ for all $1\le i\le m$, and $Bv_i=\lambda_1v_i$ for all $m+1\le i\le n$. In particular, $B$ agrees with $A$ on $U$, and so $U$ is invariant under the action of $B$.

  We claim that $p(B)=0_n$. Indeed, we will show that $p(B)v_i=0$ for all $1\le i\le n$. For $1\le i\le m$, we have $p(B)v_i=p(A)v_i=0$, since $A$ and $B$ agree on $U$ and since $U$ is invariant under this action. For $m+1\le i\le n$ we have $p(B)v_i=0$ since already $(B-\lambda_1\Id_n)v_i=0$. This shows that $p(B)=0_n$. Additionally, $A$ and $B$ agree on $m$ vectors of the above basis, hence
  \[
    \Rank(A-B) \le n-m \le \eps n,
  \]
  so $B$ is an $\eps$-approximation of $A$.
\end{example}

\begin{example}
The algebra $\Alg=F\left<x,y\right>/\left<xy\right>$ is rank-stable. 

Indeed, suppose that $A_1,A_2\in \M[n](F)$ are given and $\rank(A_1 A_2)<\eps$. Write $M=A_1A_2$, then $\rank(M)<\eps$ and $\ker A_2 \subseteq \ker M$.
Fix a basis to $\Image A_2$, say, $A_2v_1,\dots,A_2v_m$. Extend it to a basis of $F^n$ by $u_{m+1},\dots,u_{n}$. Define $D\in \M[n](F)$ by $DA_2v_i = -Mv_i$ for each $1\le i\le m$ and $Du_i=0$ for each $m+1\le i\le n$. For every $w\in F^n$ write $A_2w=\sum_{i=1}^{m} \alpha_i A_2v_i$, so $w-\sum_{i=1}^{m} \alpha_i v_i\in \ker A_2\subseteq \ker M$. Thus: \[
  DA_2w=\sum_{i=1}^{m} \alpha_i DA_2v_i=-M\left(\sum_{i=1}^{m} \alpha_i v_i\right)=-Mw.
\]
Notice that $\rank(D)\le \rank(M)<\eps$. In addition
\[
  (A_1+D)A_2=A_1A_2+DA_2=M+DA_2=0_n.
\]
Hence $B_1\coloneqq A_1+D$ and $B_2\coloneqq A_2$ satisfy $B_1B_2=0_n$, thus $\eps$-approximate $A_1,A_2$.
\end{example}

\section{Compression Lemmas} \label{sec:compression}
In this section we develop a useful ``compression machinery'' that allows us to effectively control the dimensions of exact solutions of equations by means of approximate solutions.

\begin{lemma}\label{lem:pres-linalg1}
  Let $A_1,\dots,A_d\in\M[n](F)$ and $B_1,\dots,B_d\in\M[n'](F)$ be matrices. Assume that $n'>(1+\eps d)n$ and $\Rank(\widehat{A_i}-\widehat{B_i})<\eps n$ for all $1\leq i\leq d$ and for some $\eps>0$. Then there exists an invertible matrix $E\in\GL[n'](F)$ such that:
  \begin{equation}\label{eq:many-E-hats}
      \widehat{E}\widehat{A_i} = \widehat{E^{-1}}\widehat{A_i} = \widehat{A_i}\widehat{E} = \widehat{A_i}\widehat{E^{-1}} = \widehat{A_i}
  \end{equation}
  and the last $ n'-\lfloor (1+\eps d)n \rfloor$ columns of $E^{-1}B_iE$ are zero for all $i$.
\end{lemma}

\begin{proof}
  Denote by $\set{e_1,\dots,e_{n'}}$ the standard basis of $F^{n'}$, and let $W$ be the subspace $\Span \set{e_{n+1},\dots,e_{n'}}$. 
  We construct a sequence of subspaces along with compatible ordered bases for them as follows. 
  First take $S_0=\set{e_1,\dots,e_n}$, and let $U_0=\Span(S_0)$. For each $1\le i\le d$, assume that $U_{i-1}$ has been defined with a given ordered basis $S_{i-1}$. Choose a basis $\set{v_{i,1},\dots, v_{i,k_i}}$ of: \[\Image\left(B_i\big|_{\bigcap_{j=1}^{i-1} \ker(B_j)\cap W}\right)\] and for each $v_{i,j}$ pick $w_{i,j}\in \bigcap_{j=1}^{i-1} \ker(B_j)\cap W$ such that $B_iw_{i,j}=v_{i,j}$. Now let
  \[U_i \coloneqq U_{i-1}+\Span\set{w_{i,1},\dots,w_{i,k_i}}\]
  and let $S_i$ be the ordered basis of $U_i$ obtained by appending a suitable subset of the set $\set{w_{i,1},\dots,w_{i,k_i}}$ to $S_{i-1}$.
  We thus obtain a sequence of subspaces:
  \[
    U_0=\Span\set{e_1,\dots,e_n}\sub U_1\sub\cdots\sub U_d\sub F^{n'}
  \] with compatible ordered bases $S_0\subseteq \cdots \subseteq S_d$ (namely, each $S_i$ is a prefix of $S_{i+1}$).  
  We claim that $k_i\leq \eps n$. Indeed, \[k_i\leq \Rank(B_i|_W)=\Rank\big((\widehat{B_i}-\widehat{A_i})|_W\big)\leq \eps n\] (the equality follows since $\widehat{A_i}|_W=0$). Consequently, 
  \begin{equation} \label{eq:k_i}
   \dim_F U_i-\dim_F U_{i-1}\leq k_i \le \eps n.
  \end{equation}
  Next, we claim that 
  \begin{equation} \label{eq:LHS}
  U_d+\Bigg(\bigcap_{j=1}^d\ker(B_j)\cap W\Bigg)=F^{n'}.   
  \end{equation}  
  Indeed, let $\xi\in W$. We inductively construct a sequence $\set{\xi_i}_{i=0}^d$ where $\xi_i$ belongs to $\bigcap_{j=1}^i\ker(B_j)\cap W$ as follows. First, let $\xi_0=\xi$. For each $1\le i\le d$, write $B_i\xi_{i-1}=\sum_{j=1}^{k_i}\alpha_{i,j}v_{i,j}$ for some $\alpha_{i,j}\in F$. Then define $\xi_i=\xi_{i-1}-\sum_{j=1}^{k_i}\alpha_{i,j}w_{i,j}$. By construction, for each $i$ we have $\xi_{i-1}-\xi_i\in U_i\sub U_d$ and $\xi_i\in\bigcap_{j=1}^i\ker(B_j)\cap W$. Write:
  \[
    \xi=(\xi_0-\xi_1)+\cdots+(\xi_{d-1}-\xi_d)+\xi_d,
  \]
  where $(\xi_0-\xi_1)+\cdots+(\xi_{d-1}-\xi_d)\in U_d$ and $\xi_d\in\bigcap_{j=1}^d\ker(B_j)\cap W$. This proves~\eqref{eq:LHS}.

Now, notice that: 
\begin{align*}
 m\coloneqq|S_d| & =  |S_0|+\left(|S_1|-|S_0|\right)+\cdots+\left(|S_d|-|S_{d-1}|\right) \\
 & =  n + \sum_{i=1}^{d} (\dim_F U_i - \dim_F U_{i-1})  \leq (1+\eps d) n.   
\end{align*}
by \eqref{eq:k_i}.
Complete $S_d$ to a basis $S'$ of $F^{n'}$ by appending suitable elements $b_{m+1},\dots,b_{n'}\in \bigcap_{j=1}^d\ker(B_j)\cap W$, which is indeed possible by \eqref{eq:LHS}. Let $E$ be the matrix whose columns are the vectors of $S'$. 
We claim that the last $n'-m\ge n'-(1+\eps d)n$ columns of $E^{-1}B_iE$ are zero. Indeed, for each $m+1\leq k\leq n'$ we have 
 $E^{-1}B_iEe_k=E^{-1}B_ib_k=0$. This proves that the last $(n'-m)$ columns of $E^{-1}B_iE$ are zero.

 Finally, we need to prove \eqref{eq:many-E-hats}. 
 Indeed, consider the direct sum decomposition $F^{n'}=U_0\oplus W$. On $U_0$, we have $\widehat{E}|_{U_0}=\Id_{U_0}$ and $\widehat{A_i}(U_0)\sub U_0$, so the equalities hold on $U_0$. On $W$, we have $\widehat{A_i}|_W=0$ and $\widehat{E}(W)\sub W$, so the desired equalities also holds on $W$ (since it vanishes on both sides).
\end{proof}

\begin{lemma}\label{lem:small-approx}
  Let $A_1,\dots,A_d\in\M[n](F)$ be matrices such that:
  \[
    \rank(P_j(A_1,\dots,A_d))<\delta
  \]
  for all $1\le j\le r$ and for some $\delta>0$, and let $B_1,\dots,B_d\in\M[n'](F)$ be an $\eps$-approximation of $A_1,\dots,A_d$ such that $P_j(B_1,\dots,B_d)=0$ for all $1\leq j\leq r$. Then there exists an $\eps$-approximation $C_1,\dots,C_d\in\M[m](F)$ of $A_1,\dots,A_d$ such that $P_j(C_1,\dots,C_d)=0$ for all $1\leq j\leq r$ and: \[
  \frac{1-\delta}{1+\eps d}n\le m\le (1+\eps d)n.\] Furthermore, if $\delta<(\eps d)^2$, then $(1-\eps d)n\le m\le (1+\eps d)n$.
\end{lemma}

\begin{proof}
  If $\frac{1-\delta}{1+\eps d}n\le n'\le (1+\eps d)n$, we are done. We divide the rest of the proof into two cases.

  CASE I: $n'>(1+\eps d)n$. Apply \Lref{lem:pres-linalg1} to get a matrix $E\in\GL[n'](F)$ such that $\widehat{E}\widehat{A_i}\widehat{E^{-1}}=\widehat{A_i}$, and such that the last $n'-\lfloor (1+\eps d)n \rfloor$ columns of $E^{-1}B_iE$ are zero. Write:
  \[
    E^{-1}B_iE=\begin{pmatrix}C_i & 0 \\ * & 0\end{pmatrix}
  \]
  for some $C_i\in\M[\lfloor (1+\eps d)n\rfloor](F)$. Note that:
  \[
    0=E^{-1}P_j(B_1,\dots,B_d)E=\begin{pmatrix}P_j(C_1,\dots,C_d) & 0 \\ * & P_j(0,\dots,0)\end{pmatrix}
  \]
  for all $j$, hence $P_j(C_1,\dots,C_d)=0$ for all $j$. In addition, writing $\widetilde{A_i}\in\M[\lfloor (1+\eps d)n\rfloor](F)$ to be the matrix $A_i$ padded by zeros of the same size as $C_i$, we obtain:
  \[
    \widehat{A_i}-\widehat{B_i}=\widehat{E}\widehat{A_i}\widehat{E^{-1}}-\widehat{B_i}=\widehat{E}(\widehat{A_i}-\widehat{E^{-1}}\widehat{B_i}\widehat{E})\widehat{E^{-1}}=\widehat{E}\begin{pmatrix}\widetilde{A_i}-C_i&0\\ *&0\end{pmatrix}\widehat{E^{-1}},
  \]
  so $\Rank(\widehat{A_i}-\widehat{C_i})\leq \Rank(\widehat{A_i}-\widehat{B_i})<\eps n$, proving the first case.

  CASE II: $n'<\frac{1-\delta}{1+\eps d}n$. If $P_j(0,\dots,0)=0$ for all $j$, we may define $C_i=B_i\oplus 0_{n-n'}$ and the lemma follows. Otherwise, note that $(1+\eps d)n'<(1-\delta)n<n$, so we may apply \Lref{lem:pres-linalg1} with reversed roles to get a matrix $E$ such that the last $n-\lfloor (1+\eps d)n' \rfloor$ columns of $E^{-1}A_iE$ are zero. Write:
  \[
    E^{-1}A_iE=\begin{pmatrix}A_i' & 0 \\ * & 0\end{pmatrix}
  \]
  for some $A_i'\in\M[\lfloor(1+\eps d)n'\rfloor](F)$. Note that:
  \[
    E^{-1}P_j(A_1,\dots,A_d)E=\begin{pmatrix}P_j(A_1',\dots,A_d') & 0 \\ * & P_j(0,\dots,0)\end{pmatrix}.
  \]
  By assumption, $P_j(0,\dots,0) = \theta_j\cdot \Id_{n-\lfloor(1+\eps d)n'\rfloor}\ne 0$ for some $j$ and $\theta_j\in F^\times$, is a non-zero scalar matrix. Hence:
  \begin{align*}
    \frac{n-\lfloor(1+\eps d)n'\rfloor}{n} = \frac{1}{n}\Rank\left(\theta_j\cdot \Id_{n-\lfloor(1+\eps d)n'\rfloor}\right) \le \rank(P_j(A_1,\dots,A_d))<\delta,
  \end{align*}
  a contradiction, since $n-\lfloor(1+\eps d)n'\rfloor\geq n-(1+\eps d)n'>\delta n$.\medskip

  Finally, note that if $\delta<(\eps d)^2$, then
  \[
    \frac{1-\delta}{1+\eps d}n\ge\frac{1-(\eps d)^2}{1+\eps d}n=(1-\eps d)n,
  \]
  proving the last statement of the lemma.
\end{proof}

The final ingredient we need is an estimation on the rank differences when substituting matrices into a polynomial. Recall that for a polynomial $f(x_1,\dots,x_d)$ and matrices $A_1,\dots,A_d\in\M[n](F)$, we denote by $\widehat{f}(A_1,\dots,A_d)\in\M[\N](F)$ the matrix obtained by applying $\,\widehat{\cdot}\,$ to $f(A_1,\dots,A_d)$.

\begin{lemma}\label{lem:computing polyrank}
  Let $A_1,\dots,A_d\in\M[n](F)$ and $B_1,\dots,B_d\in\M[n'](F)$ be matrices. Suppose that $\Rank(\widehat{A_i}-\widehat{B_i})<\lambda$ for all $1\leq i\leq d$. Let
  \[
    f(x_1,\dots,x_d)=\sum_{j=1}^l c_j x_{\alpha_{j,1}}\cdots x_{\alpha_{j,p_j}}
  \]
  be a polynomial, and let $m$ be such that $p_j\leq m$ for all $1\leq j\leq l$. Then:
  \[
    \Rank(\widehat{f}(A_1,\dots,A_d)-\widehat{f}(B_1,\dots,B_d))<lm\lambda+\left|n-n'\right|.
  \]
  If $f$ has a zero constant term, the summand $\left|n-n'\right|$ can be omitted.
\end{lemma}

\begin{proof}
  Suppose that $\Rank(C-D)<t$ for some matrices $C,D\in\M[\N](F)$. Then:
  \begin{align*}
\Rank(\widehat{A_i}C-\widehat{B_i}D) & = \Rank((\widehat{A_i}-\widehat{B_i})C+\widehat{B_i}(C-D)) \\ & \leq \Rank(\widehat{A_i}-\widehat{B_i})+\Rank(C-D)<\lambda+t.
  \end{align*}
  By induction, for each non-empty monomial $x_{\alpha_{j,1}}\cdots x_{\alpha_{j,p_j}}$ we have that:
  \[
    \Rank(\widehat{A_{\alpha_{j,1}}}\cdots \widehat{A_{\alpha_{j,p_j}}} - \widehat{B_{\alpha_{j,1}}}\cdots \widehat{B_{\alpha_{j,p_j}}}) < p_j\lambda \leq m\lambda.
  \]
  Note that a scalar monomial is interpreted as a scalar matrix of the matching size. If $n=n'$, it cancels out in the difference $\widehat{f}(A_1,\dots,A_d)-\widehat{f}(B_1,\dots,B_d)$. Otherwise, the difference between the constant monomials in $\widehat{f}(A_1,\dots,A_d)$ and $\widehat{f}(B_1,\dots,B_d)$ is a diagonal matrix with $\left|n-n'\right|$ non-zero entries, which has rank $\left|n-n'\right|$.

  Summing up,
  \begin{align*}
     & \Rank(\widehat{f}(A_1,\dots,A_d) - \widehat{f}(B_1,\dots,B_d)) = \\  & \Rank\Bigg(\sum_{j=1}^l c_j A_{\alpha_{j,1}}\cdots A_{\alpha_{j,p_j}}-\sum_{j=1}^l c_j B_{\alpha_{j,1}}\cdots B_{\alpha_{j,p_j}}\Bigg) \leq \\
     & \sum_{j=1}^l \Rank\left(A_{\alpha_{j,1}}\cdots A_{\alpha_{j,p_j}} - B_{\alpha_{j,1}}\cdots B_{\alpha_{j,p_j}} \right) \leq lm\lambda+\left|n-n'\right|,
  \end{align*}
  as claimed.
\end{proof}

\section{Independence of Presentation} \label{sec:independence}

The following is natural but surprisingly non-trivial, and more involved than the group case.

\begin{theorem}\label{thm:two-pres}
  Rank-stability of a finitely presented algebra does not depend on its presentation.
\end{theorem}

\begin{proof}
  Let $\Alg$ be a finitely presented algebra with two presentations:
  \[
    F\left<x_1,\dots,x_d\right>/\left<P_1,\dots,P_r\right>\cong \Alg\cong F\left<y_1,\dots,y_e\right>/\left<Q_1,\dots,Q_s\right>,
  \]
  and assume that the left hand side presentation is rank-stable.

  The idea is as follows. If we have an approximate solution to $Q_1,\dots,Q_s$, we may use an algebra isomorphism to translate it to an approximate solution to $P_1,\dots,P_r$. The rank-stability of the first presentation lets us find an approximation of the new matrices. We then translate it back to an approximation of the original matrices by the inverse isomorphism.

  For each $1\le i\le d$ let $G_i$ be an $e$-variate non-commutative polynomial such that $x_i\mapsto G_i(\vec{y})$ defines an isomorphism. Similarly, for each $1\le  i\le e$ let $F_i$ be a $d$-variate non-commutative polynomial such that $y_i\mapsto F_i(\vec{x})$ defines the inverse isomorphism. Let $S$ bound from above the product of the maximal length of any monomial appearing in any $F_i$, times the maximal number of monomials appearing in any $F_i$.

  Since $F_i$ and $G_i$ define inverse isomorphisms between the two presentations of~$\Alg$, for each $1\le j\le r$ the equations
  \begin{align}
    F_j(G_1(\vec{y}),\dots,G_d(\vec{y})) & = y_j + \sum_{l=1}^{K_j} a_{j,l} Q_{m_{j,l}}(\vec{y}) b_{j,l} \label{eq:srs-pres0}\\
    P_j(G_1(\vec{y}),\dots,G_d(\vec{y})) & = \sum_{l=1}^{M_j} u_{j,l} Q_{\mu_{j,l}}(\vec{y}) v_{j,l} \label{eq:srs-pres1}
  \end{align}
  hold in the free algebra $F\left<y_1,\dots,y_e\right>$ for some $a_{j,l},b_{j,l},u_{j,l},v_{j,l}\in F\left<y_1,\dots,y_e\right>$. Similarly, for each $1\le j\le s$ the equation
  \begin{equation}\label{eq:srs-pres2}
      Q_j(F_1(\vec{x}),\dots,F_e(\vec{x})) = \sum_{l=1}^{N_j} w_{j,l} P_{\nu_{j,l}}(\vec{x}) z_{j,l}
  \end{equation}
  holds in the free algebra $F\left<x_1,\dots,x_d\right>$ for some $w_{j,l},z_{j,l}\in F\left<x_1,\dots,x_d\right>$. Let us assume that all $M_j,N_j$'s are bounded from above by some integer $T$, and that the $K_j$'s are bounded by $R$.

  Let $\eps>0$ be given. Let $\delta>0$ be an appropriate $\delta$ from \Dref{def:rank-stab} of rank-stability, according to the presentation $F\left<x_1,\dots,x_d\right>/\left<P_1,\dots,P_r\right>$ and with respect to $\eps'\coloneqq\frac{\eps}{2(S+d)}$. We also assume that $\delta<(\eps' d)^2$.

  Let $\vec{A}=(A_1,\dots,A_e)\in\M[n](F)$ be matrices satisfying $\rank(Q_i(\vec{A}))<\min\set{\frac{\delta}{T},\frac{\eps}{2R}}$ for any $1\le i\le s$. It follows from \eqref{eq:srs-pres1} that for any $1\le j\le r$ we have
  \[
    \rank \left(P_j(G_1(\vec{A}),\dots,G_d(\vec{A}))\right) < T\cdot \delta/T=\delta.
  \]
  By the rank-stability of $F\left<x_1,\dots,x_d\right>/\left<P_1,\dots,P_r\right>$, there exist matrices $B_1,\dots,B_d\in\M[n'](F)$ such that $P_j(B_1,\dots,B_d)=0$ for any $1\le j\le r$, and $\Rank(\widehat{G_i}(\vec{A})-\widehat{B_i})<\eps' n$ for any $1\le i\le d$. By \Lref{lem:small-approx}, as $\delta<(\eps' d)^2$, we may assume that $(1-\eps'd)n\le n'\le (1+\eps' d)n$, so $\left|n'-n\right|\le \eps' dn$.
  
  Let $C_i=F_i(\vec{B})$ for each $1\le i\le e$. We will show that $\vec{C}=(C_1,\dots,C_e)$ is a solution to $Q_1,\dots,Q_s$ that $\eps$-approximates $\vec{A}$. Indeed, notice that $Q_j(\vec{C})=0$ for each $1\le j\le s$ by \eqref{eq:srs-pres2}. We claim that $\vec{C}$ is an $\eps$-approximation of $\vec{A}$, since
  \begin{align*}
    \Rank\left(\widehat{C_i} - \widehat{A_i}\right) & \le \Rank\left(\widehat{C_i} - \widehat{F_i}(G_1(\vec{A}),\dots,G_d(\vec{A}))\right) \\
      &\qquad + \Rank\left(\widehat{F_i}(G_1(\vec{A}),\dots,G_d(\vec{A})) - \widehat{A_i}\right) \\
    & = \Rank\left(\widehat{F_i}(\vec{B}) - \widehat{F_i}(G_1(\vec{A}),\dots,G_d(\vec{A}))\right) \\
    &\qquad + \Rank\left(F_i(G_1(\vec{A}),\dots,G_d(\vec{A})) - A_i\right).
  \end{align*}
  The first term is bounded by $S\eps' n+\left|n'-n\right|\le S\eps' n+\eps' dn=\frac{\eps}{2}n$ by \Lref{lem:computing polyrank}, and the second term is bounded by $R\cdot\frac{\eps}{2R}n=\frac{\eps}{2}n$ by \eqref{eq:srs-pres0}. Therefore $\Rank\big(\widehat{C_i}-\widehat{A_i}\big)\le \eps n$, concluding the proof.
\end{proof}

\section{Finite-Dimensional Algebras} \label{sec:finite-dimensional}

Our goal in this section is to prove the following:

\begin{theorem}\label{thm:fin-dim}
  Let $\Alg$ be a finite-dimensional algebra over a field $F$. Then $\Alg$ is rank-stable.
\end{theorem}

Our proof strategy is partially analogous to the proof that finite groups are P-stable \cite{GlebskyRivera09}. 
Recall that every finite-dimensional algebra over a field is finitely presented.

Let $\Alg$ be a finite-dimensional algebra. Fix an $F$-vector space basis $\{1,a_1,\dots,a_d\}$ of $\Alg$. For each $1\leq i,j\leq d$ we can write $a_i a_j = \sum_{k=1}^d\alpha_{ijk} a_d + \beta_{ij}$ for some scalars $\alpha_{ijk},\beta_{ij}\in F$.
Then $\Alg\cong F\left<x_1,\dots,x_d\right>/\left<P_{11},P_{12},\dots,P_{dd}\right>$, where $P_{ij}$ is defined by
\[
    P_{ij} = x_ix_j - \sum_{k=1}^d\alpha_{ijk}x_d - \beta_{ij}.
\]
Indeed, let $\pi\colon F\left<x_1,\dots,x_d\right>\twoheadrightarrow \Alg$ by $\pi(x_1)=a_1,\dots,\pi(x_d)=a_d$. Clearly, each $P_{ij}$ belongs in $\ker \pi$, so $\pi$ induces a surjection $\overline{\pi}\colon F\left<x_1,\dots,x_d\right>/\left<P_{11},P_{12},\dots,P_{dd}\right>\twoheadrightarrow \Alg$. Moreover, notice that $F\left<x_1,\dots,x_d\right>/\left<P_{11},P_{12},\dots,P_{dd}\right>$ is spanned by (the cosets of) $1,x_1,\dots,x_d$, and since $\dim_F \Alg = d+1$ then $\overline{\pi}$ must be an isomorphism.

Fix $\eps>0$. Let $\delta>0$ be a constant that will be specified later. Let $A_1,\dots,A_d\in\M[n](F)$ be matrices such that
\[
    \rank(P_{ij}(A_1,\dots,A_d))<\delta
\]
for every $1\le i,j\le d$.

\begin{notation}
      For a matrix $B\in\M[n](F)$ and a subspace $U\le F^n$, we denote the pre-image of $U$ under $B$ by $B^{-1}U=\set{v\in F^n\mid Bv\in U}$.
\end{notation}

Notice that $\dim_F B^{-1}U\ge\dim_F U$. Indeed, pick a basis $Bv_1,\dots,Bv_r$ to $U\cap \Image(B)$ and complete it to a basis of $U$ by some $u_1,\dots,u_s$. Then $\dim_F \ker(B)\geq s$ and $\ker(B)\oplus\Span \{v_1,\dots,v_r\}\subseteq B^{-1}U$.

\begin{fact}
Recall that if $U_1,\dots,U_d\leq F^n$, then $\dim_F \bigcap_{i=1}^d U_i \geq  n - \sum_{i=1}^d (n-\dim_F U_i)$. This follows since the natural linear map $F^n/\bigcap_{i=1}^d U_i \hookrightarrow F^n/U_1 \oplus \cdots \oplus F^n/U_d$ is injective.
\end{fact}

Let $U=\bigcap_{i,j=1}^d \ker P_{ij}(A_1,\dots,A_d)$, and notice that $\dim_F U\ge (1-d^2\delta)n$. Let
\[
W = U\cap \bigcap_{i=1}^d A_i^{-1}U,
\]
and notice that $W\subseteq U$.

\begin{lemma}\label{lem:fin-dim-act}
There is a well-defined action $\Alg\curvearrowright W$ via $x_iw=A_iw$.
\end{lemma}

\begin{proof}
First we show that $A_iW\subseteq W$ for each $1\leq i\leq d$. Let $w\in W$. By the definition of $W$, we have $A_iw\in U$, so it remains to show that $A_iw\in A_j^{-1}U$, i.e., $A_jA_iw\in U$, for all $1\le j\le d$. Recall that $w\in U$, hence $P_{ji}(A_1,\dots,A_d)w=0$, meaning
\[
    A_jA_iw = \sum_{k=1}^d\alpha_{jik}A_kw+\beta_{ji}w.
\]
Each $A_kw\in U$ since $w\in W\sub A_k^{-1}U$, and $w\in U$, hence $A_jA_iw\in U$.

Therefore $F\left<x_1,\dots,x_d\right>\curvearrowright W$ via $x_iw=A_iw$. Furthermore, since $W\subseteq U$, it follows that $P_{11},P_{12},\dots,P_{dd}$ annihilate $W$. Hence this action reduces to an action $\Alg\curvearrowright W$.
\end{proof}

\begin{lemma}\label{lem:fin-dim-good}
  There exists a constant $c>0$, which depends only on $d$ but not on $\delta$ or $\vec{A}$, such that $\dim_F W\ge (1-c\delta)n$.
\end{lemma}

\begin{proof}
We have
\begin{align*}
    \dim_F W & \geq \dim_F U - \sum_{i=1}^d(n-\dim_F A_i^{-1}U) \geq \dim_F U - \sum_{i=1}^d(n-\dim_F U) \\
    & = (d+1)\dim_F U - dn \ge (d+1)(1-d^2\delta)n - dn = (1-(d+1)d^2\delta)n
\end{align*}
as required.
\end{proof}

We are now ready to conclude the proof of \Tref{thm:fin-dim}. Since $\Alg$ is finite-dimensional, it has a faithful representation $\varphi\colon \Alg\hookrightarrow\M[s](F)$ for some constant $s\ge 1$. In fact, we do not use the faithfulness of this representation.
Let $C_i=\varphi(x_i)$ for each $1\leq i\leq d$, so $P_{ij}(C_1,\dots,C_d)=0_s$ for each $1\leq i,j\leq d$.

Let $\eps>0$. Take $0<\delta<\min\set{\eps/2c,\eps/2s}$ for $c$ from \Lref{lem:fin-dim-good}, which is possible since $c$ and $s$ do not depend on $\delta$. This choice of $\delta$ guarantees that $k \coloneqq \dim_F W\ge\allowbreak (1-\eps/2)n$. In addition, note that if $n<2s/\eps$, then $\rank(P_{ij}(A_1,\dots,A_d))<\eps/2s<1/n$, so $P_{ij}(A_1,\dots,A_d)=0_n$ for all $i,j$. Therefore $\vec{A}=(A_1,\dots,A_d)$ is an exact solution to $P_{11},P_{12},\dots,P_{dd}$. We henceforth assume that $n\ge 2s/\eps$.

Let $S=\set{v_1,\dots,v_k}$ be a basis of $W$. Complete $S$ to a basis $S'=\{v_1,\dots,v_k,\allowbreak v_{k+1},\dots,v_n\}$ of $F^n$, and place the vectors of $S'$ as the columns of a matrix $E$.

For each $1\leq i\leq d$, notice that $E^{-1}A_iE$ is a block upper triangular matrix, since $W=\Span_F \{Ee_1,\dots,Ee_k\}$ is $A_i$-invariant. Denote by $B_i\in\M[k](F)$ the submatrix of $E^{-1}A_iE$ obtained by taking the first $k$ rows and columns of $E^{-1}A_iE$. 
Define a block diagonal matrix:
\[
    B_i'=\begin{pmatrix}
        B_i \\ & C_i \\ & & \ddots \\ & & & C_i
    \end{pmatrix}\in\M[n'](F)
\]
for $n'$ the smallest integer such that $n'\geq n$ and $n'\equiv k\pmod{s}$. Note that $\widehat{E^{-1}A_iE}$ and $\widehat{B_i'}$ act the same way on $e_1,\dots,e_k$, so:
\[
    \Rank(\widehat{E^{-1}A_iE}-\widehat{B_i'})\le n'-k = (n-k)+(n'-n)\le \eps n/2+s\le\eps n.
\]

We also claim that for each $1\leq i,j\leq d$ we have $P_{ij}(B_1',\dots,B_d')=0_{n'}$. Indeed, let us show that $P_{ij}(B_1',\dots,B_d')e_l=0$ for every $1\leq l\leq n'$. If $1\le l\le k$, then:
\[
    P_{ij}(B_1',\dots,B_d')e_l=E^{-1}P_{ij}(A_1,\dots,A_d)Ee_l=E^{-1}P_{ij}(A_1,\dots,A_d)v_l=0
\]
since $v_l\in W\subseteq U$. If $k+1\le l\le n'$, then:
\[
    P_{ij}(B_1',\dots,B_d')e_l=(0,\dots,0,P_{ij}(C_1,\dots,C_d)e_p,0\dots,0)=0_{n'}
\]
for some $1\leq p\leq s$ since $P_{ij}(C_1,\dots,C_d)=0_s$. Specifically, $p=\left((l-k-1)\mod{s}\right)+1$.

Let $E'=E\oplus \Id_{n'-n}$. Note that $\widehat{E'}\widehat{E^{-1}A_iE}\widehat{(E')^{-1}}=\widehat{A_i}$ and $\widehat{E'}\widehat{B_i'}\widehat{(E')^{-1}}=\stackon[-8pt]{\ensuremath{E'B_i'(E')^{-1}}}{\vstretch{1.5}{\hstretch{9.0}{\widehat{\phantom{\;}}}}}$. Therefore,
\begin{align*}
\Rank(\widehat{A_i}-\stackon[-8pt]{\ensuremath{E'B_i'(E')^{-1}}}{\vstretch{1.5}{\hstretch{9.0}{\widehat{\phantom{\;}}}}}) & =\Rank(\widehat{E'}\widehat{E^{-1}A_iE}\widehat{(E')^{-1}}-\widehat{E'}\widehat{B_i'}\widehat{(E')^{-1}})\\
    &\le\Rank(\widehat{E^{-1}A_iE}-\widehat{B_i'})\le\eps n.
\end{align*}
The matrices $E'B_1'(E')^{-1},\dots,E'B_d'(E')^{-1}$ are therefore a solution to $P_{11},P_{12},\dots,P_{dd}$ and $\eps$-approximate $A_1,\dots,A_d$.

\section{Algebra Stability and Group Stability} \label{sec:group}

In this section we relate the rank-stability of a group with its group algebra. We recall the definition from Elek and Grabowski~\cite{ElekGrabowski21} (studied therein for unitary matrices).

\begin{definition}
  Let $G=\left<x_1,\dots,x_d\mid P_1,\dots,P_r\right>$ be a finitely presented group, and let $F$ be a field. We say that $G$ is \textbf{rank-stable over $F$} if for every $\eps>0$ there exists $\delta>0$ such that the following holds: for any $n\ge 1$ and $A_1,\dots,A_d\in\GL[d](F)$ such that $\rank (P_j(A_1,\dots,A_d)-\Id_n)<\delta$ for all $1\le j\le r$, there exist $n'\ge 1$ and $B_1,\dots,B_d\in\GL[n'](F)$ such that $\Rank(\widehat{A_i}-\widehat{B_i})<\eps n$ for all $1\le i\le d$, and $P_j(B_1,\dots,B_d)=\Id_{n'}$ for all $1\le j\le r$.
\end{definition}

\begin{lemma}\label{lem:inv-diff}
    Let $A\in\GL[n](F)$ and $B\in\GL[n'](F)$ be invertible matrices. Then
    \[
        \Rank(\widehat{A^{-1}} - \widehat{B^{-1}}) \le \Rank(\widehat{A} - \widehat{B}) + 2\left|n'-n\right|.
    \]
\end{lemma}

\begin{proof}
    Without loss of generality, we assume $n'\ge n$. Let $A'=A\oplus\Id_{n'-n}\in\GL[n'](F)$. First note that
    \begin{align*}
        \Rank(\widehat{A^{-1}} - \widehat{B^{-1}}) &\le \Rank(\widehat{A^{-1}} - \widehat{(A')^{-1}}) + \Rank(\widehat{(A')^{-1}} - \widehat{B^{-1}}) \\
        &= \Rank((A')^{-1} - B^{-1}) + (n'-n).
    \end{align*}
    Next, since $A'$ and $B$ are invertible and of the same dimensions,
    \begin{align*}
        \Rank((A')^{-1} - B^{-1}) &= \Rank(A'((A')^{-1} - B^{-1})B) = \Rank(B-A') = \Rank(\widehat{B}-\widehat{A'}) \\
        & \le \Rank(\widehat{B}-\widehat{A}) + \Rank(\widehat{A}-\widehat{A'}) = \Rank(\widehat{A}-\widehat{B}) + (n'-n)
    \end{align*}
    which proves the lemma.
\end{proof}

\begin{theorem}\label{thm:group-alg}
  Let $G$ be a finitely presented group, and let $F$ be a field. Then $G$ is rank-stable over $F$ if and only if its group algebra $F[G]$ is rank-stable.
\end{theorem}

\begin{proof}
  Fix a presentation $G=\left<x_1,\dots,x_d\mid P_1,\dots,P_r\right>$. Then $F[G]$ has the presentation
  \[
    F[G]=F\left<x_1,\dots,x_d,y_1,\dots,y_d\right>/\left<Q_i(\vec{x},\vec{y})-1,x_jy_j-1,y_jx_j-1\right>_{1\le i\le r,1\le j\le d},
  \]
  where $Q_i$ is obtained from $P_i$ by replacing each $x_j^{-1}$ by $y_j$. Let $C$ denote twice the maximal degree of the $Q_i$'s.
  
  Suppose first that $G$ is rank-stable over $F$. Let $\eps>0$, and take $\delta$ according to the rank-stability of $G$ with $\eps'=\frac{\eps}{2(1+2d)}$, and assume that $\delta<(\eps'd)^2$. Let $A_1,\dots,A_d,A_1',\dots,A_d'\in\M[n](F)$ be matrices that are a $\delta'$-approximate solution of the relations of $F[G]$, where $\delta'=\min\set{\frac{\delta}{C},\frac{\eps}{10}}$. In particular, for any $j$ we have
  \[
    1=\rank(\Id_n)\le\rank(A_jA_j')+\rank(A_jA_j'-\Id_n)\le\min\set{\rank(A_j),\rank(A_j')}+\delta',
  \]
  so $\rank(A_j),\rank(A_j')\ge 1-\delta'$. Therefore one may write $A_j=U_j+E_j$ and $A_j'=U_j'+E_j'$, where each $U_j$ and $U_j'$ is invertible and $\rank(E_j),\rank(E_j')\le\delta'$. Using $\rank(A_jA_j'-\Id_n)\le\delta'$ again, it follows that
  \begin{align*}
      \rank(U_j'-U_j^{-1}) &= \rank(U_jU_j'-\Id_n) = \rank((A_j-E_j)(A_j'-E_j')-\Id_n) \\
      & \le \rank(A_jA_j'-\Id_n) + \rank(A_jE_j') + \rank(E_jA_j') + \rank(E_jE_j')\le 4\delta'.
  \end{align*}
Recall also that for any $j$ we have $\rank(Q_j(\vec{A},\vec{A'}))\le\delta'$. Therefore \Lref{lem:computing polyrank} shows that
  \[
    \rank(P_j(\vec{U})-\Id_n)=\rank(Q_j(\vec{U},\vec{U^{-1}})-\Id_n)\le C\delta'\le\delta.
  \]
  By the rank-stability of~$G$, there exist $n'\ge 1$ and invertible matrices $V_1,\dots,V_d\in\GL[n'](F)$ such that $P_j(\vec{V})=\Id_{n'}$ for all $j$, and such that $\Rank(\widehat{V_j}-\widehat{U_j})\le\eps'n$. Further, by \Lref{lem:small-approx}, we may assume $\left|n'-n\right|\le\eps'dn$.

  We show that $(\vec{V},\vec{V^{-1}})$ is an $\eps$-approximation of $(\vec{A},\vec{A'})$. Indeed,
  \[
    \Rank(\widehat{A_j}-\widehat{V_j})\le\Rank(\widehat{A_j}-\widehat{U_j})+\Rank(\widehat{U_j}-\widehat{V_j})\le\delta'n+\eps'n\le\eps n.
  \]
  By \Lref{lem:inv-diff} we have
  \[
    \Rank(\widehat{U_j^{-1}}-\widehat{V_j^{-1}}) \le \Rank(\widehat{U_j}-\widehat{V_j}) + 2\left|n'-n\right| \le (\eps'+2\eps' d)n,
  \]
  and thus
  \begin{align*}
    \Rank(\widehat{A_j'} - \widehat{V_j^{-1}}) &\le \Rank(\widehat{A_j'}-\widehat{U_j'})+\Rank(\widehat{U_j'}-\widehat{U_j^{-1}})+\Rank(\widehat{U_j^{-1}} - \widehat{V_j^{-1}}) \\
    &\le \delta'n+4\delta'n+\eps'(1+2d)n\le\eps n.
  \end{align*}
  As $Q_j(\vec{V},\vec{V^{-1}})=\Id_{n'}$ for all $j$, we proved the rank-stability of $F[G]$.\bigskip

  Suppose now that $F[G]$ is rank-stable. Let $\eps>0$, and let $\delta>0$ from the rank-stability of $F[G]$ with respect to $\eps$. Assume that $A_1,\dots,A_d\in\GL[n](F)$ are invertible matrices such that $\rank(P_j(\vec{A})-\Id_n)<\delta$ for all $j$. For each $j$ set $B_j=A_j^{-1}$; then $A_jB_j=B_jA_j=\Id_n$ and $\rank(Q_j(\vec{A},\vec{B})-\Id_n)<\delta$. By the rank-stability of $F[G]$, there are matrices $U_1,\dots,U_d,V_1,\dots,V_d\in\M[n'](F)$ such that $Q_i(\vec{U},\vec{V})=\Id_{n'}$, $U_jV_j=V_jU_j=\Id_{n'}$, and $\Rank(\widehat{A_j}-\widehat{U_j})\le\eps n$. In particular, $V_j=U_j^{-1}$, so $P_i(\vec{U})=\Id_{n'}$ as well, and $\vec{U}$ is the desired $\eps$-approximation of $\vec{A}$. This proves that~$G$ is rank-stable over~$F$.
\end{proof}

\section{Stability Under Various Constructions} \label{sec:constructions}

In this section, we study how rank-stability behaves under three natural constructions: free products, direct products, and matrix algebras. We show that, in each case, the rank-stability of all of the components is equivalent to the rank-stability of the resulting algebra.

Notice that a free product of P-stable groups is P-stable, but Ioana \cite{Ioana20} proved that the direct product is not necessarily P-stable, answering a question of Becker, Lubotzky and Thom \cite{BeckerLubotzkyThom19}. 
Ioana \cite{Ioana21} further showed that the direct product of Hilbert--Schmidt stable groups is not necessarily Hilbert--Schmidt stable. However, if at least one of the factors of that direct product has Kazhdan's property (T), then the direct product is Hilbert--Schmidt stable by a result of de la Salle \cite{delaSalle22}.
It would be interesting to find an example of finitely presented, residually finite rank-stable groups (resp.\ algebras) whose direct product (resp.\ tensor product) is not rank-stable.

We add the assumption that the algebras have a finite-dimensional representation, to avoid obstructions such as \Exref{ex:vac-stab}.

Recall that the \textbf{Kronecker product} of  matrices $A\in\M[n](F)$ and $B\in\M[n'](F)$ is the matrix $A\otimes B\in\M[nn'](F)$ defined by
\begin{equation}\label{eq:kronecker-product}
    A\otimes B=\begin{pmatrix} a_{11}B & \cdots & a_{1n}B \\ \vdots & \ddots & \vdots \\ a_{n1}B & \cdots & a_{nn}B\end{pmatrix}.
\end{equation}
In particular, the matrix $\Id_n\otimes B$ is an $n\times n$ block scalar matrix, where each diagonal block is the matrix $B$.

\subsection{Free products}

\begin{proposition} \label{prop:free_product}
Let $\Alg,\Alg'$ be finitely presented algebras with finite-dimensional representations. 
Then the free product $\Alg*\Alg'$ is rank-stable if and only if both $\Alg,\Alg'$ are rank-stable.
\end{proposition}
\begin{proof}
Let \[ \mathcal{R}=\{n\suchthat\exists \rho\colon \Alg\rightarrow \M[n](F)\},\quad \mathcal{R}'=\{n\suchthat\exists \rho\colon \Alg'\rightarrow \M[n](F)\}\subseteq \N\]
be the sets of dimensions of the finite-dimensional representations of $\Alg$ and $\Alg'$.

Then $\mathcal{R}, \mathcal{R}'$ are (additive) subsemigroups of $\mathbb{N}$, and therefore there exist $g,g'\in \N$ such that
\[ |\mathcal{R}\triangle g\N|,|\mathcal{R}'\triangle g'\N|<\infty.\]
 
Let $c=\text{gcd}(g,g')$ and fix $k,k'\in \mathbb{Z}_{\geq 0}$ such that $kg-k'g'=c$. Furthermore, we may assume that $k,k'$ are sufficiently large such that there exist representations $\rho\colon \Alg\rightarrow \M[kg](F)$ and $\rho'\colon \Alg'\rightarrow \M[k'g'](F)$. 
Fix presentations:
\[ \Alg\cong F\left<x_1,\dots,x_d\right>/\left<P_1,\dots,P_r\right>,\quad \Alg'\cong F\left<y_1,\dots,y_e\right>/\left<Q_1,\dots,Q_s\right>. \]
Then $\Alg*\Alg'\cong F\left<x_1,\dots,x_d,y_1,\dots,y_e\right>/\left<P_1,\dots,P_r,Q_1,\dots,Q_s\right>$. 
Suppose that $\Alg,\Alg'$ are rank-stable; fix $\varepsilon>0$ and matrices $A_1,\dots,A_d,A'_1,\dots,A'_e\in \M[n](F)$ such that their normalized rank satisfy $\rank P_i(\vec{A}),\rank Q_j(\vec{A'}) < \delta$ for $\delta=\min\{\delta_{\Alg},\delta_{\Alg'}\}$ arising from the rank-stability of $\Alg,\Alg'$.

By the rank-stability of $\Alg,\Alg'$ and \Lref{lem:small-approx}, there exist $B_1,\dots,B_d\in \M[m](F)$ and $B'_1,\dots,B'_e\in \M[m'](F)$ such that:
\begin{itemize}
    \item $P_i(\vec{B}),Q_j(\vec{B'})=0$ for all $1\leq i\leq r,1\leq j\leq s$;
    \item $\Rank(\widehat{B_i}-\widehat{A_i}),\Rank (\widehat{B'_j}-\widehat{A'_j})<\varepsilon n$ for all $1\leq i\leq d,1\leq j\leq e$;
    \item $\frac{1-\delta}{1+\varepsilon (d+e)}n\leq m,m'\leq (1+\varepsilon(d+e))n$.
\end{itemize}
By decreasing $\delta$ if needed, we may assume that $n$ is sufficiently large, such that $m,m'$ are also sufficiently large such that $m\in g\N,m'\in g'\N$. 
Without loss of generality, assume that $m\leq m'$; obviously, $c|m'-m$, so we can write $m'-m=qc=qkg-qk'g'$ for some $q\in \mathbb{Z}_{\geq 0}$, so $m+qkg=m'+qk'g'$. Furthermore, notice that $qk,qk'\leq (m'-m)\cdot \max\{k,k'\}$, where $k,k'$ are constants depending only on the algebras $\Alg,\Alg'$. Let
\[ C_i=B_i\oplus \left(\Id_q \otimes \rho(x_i)\right),\ C'_j=B'_j\oplus \left(\Id_q\otimes \rho'(y_j)\right) \in \M[m+qkg](F)=\M[m'+qk'g'](F) \]
for every $1\leq i\leq d,1\leq j\leq e$.
Define an $m+qkg=m'+qk'g'$-dimensional representation of $\Alg*\Alg'$ by $x_i\mapsto C_i$ and $y_j\mapsto C'_j$. Then:
\begin{align*}
\Rank \left(\widehat{C_i}-\widehat{A_i}\right) & \leq  \Rank \left(\widehat{C_i}-\widehat{B_i}\right)  + \Rank \left(\widehat{B_i}-\widehat{A_i}\right) \\
& \leq  qkg+\varepsilon n \\
& \leq  (m'-m)\cdot \max\{k,k'\}\cdot g +\varepsilon n \\
& \leq  \left((1+\varepsilon(d+e))n-\frac{1-\delta}{1+\varepsilon (d+e)}n\right)\cdot \max\{k,k'\} \cdot g + \varepsilon n \\
& \leq  \left(2\varepsilon(d+e)+\varepsilon^2(d+e)^2+\delta+\varepsilon\right) \cdot \max\{k,k'\} \cdot g \cdot n
\end{align*}
and similarly for $\Rank(\widehat{C'_i}-\widehat{A'_i})$. So, taking $\delta<\eps$ in the first place, we obtain a solution to all of $P_1,\dots,P_r,Q_1,\dots,Q_s$ which $\alpha\varepsilon$-approximates $A_1,\dots,A_d,\allowbreak A'_1,\dots,A'_e$ (where $\alpha$ is a constant depending only on $d,e,k,k',g,g'$), proving the rank-stability of $\Alg*\Alg'$.

Conversely, assume that $\Alg*\Alg'$ is rank-stable and let us show that $\Alg$ is rank-stable. Let $\varepsilon>0$ be given. Pick an arbitrary $t\in \mathcal{R}'$ and let $\psi\colon \Alg'\rightarrow M_t(F)$.
Let $\delta>0$ be associated with $\varepsilon$ from the rank-stability of $\Alg*\Alg'$ and let $\overline{\delta}\coloneqq\frac{\delta}{t}$. 
Suppose that $A_1,\dots,A_d\in \M[n](F)$ satisfy $\rank P_i(\vec{A}) <\overline{\delta}$ for each $1\leq i\leq r$. If $n\leq \frac{t}{\delta}$ then $\vec{A}$ is in fact an exact solution to $P_1,\dots,P_r$. Otherwise $n>\frac{t}{\delta}$, and write $n=qt+u$ with $0\le u<t$ and $q\in \mathbb{Z}_{\geq 0}$. Let $A'_j=\left(\Id_q \otimes \psi(y_j)\right) \oplus 0_u\in \M[n](F)$ for each $1\leq j\leq e$ and notice that $\rank Q_k(\vec{A'})\leq\frac{u}{n}<\frac{t}{n}<\delta$ for every $1\leq k\leq s$. So by the rank-stability of $\Alg*\Alg'$ there exists an exact solution $B_1,\dots,B_d,B'_1,\dots,B'_e$ to $P_1,\dots,P_r,Q_1,\dots,Q_s$ which $\varepsilon$-approximates $A_1,\dots,A_d,A'_1,\dots,A'_e$, so in particular $B_1,\dots,B_d$ approximate $A_1,\dots,A_d$ as desired.
\end{proof}

\subsection{Direct products}
Next, we turn to direct products:

\begin{theorem} \label{thm:direct-product}
Let $\Alg,\Alg'$ be finitely presented algebras with finite-dimensional representations. 
Then the direct product $\Alg\times\Alg'$ is rank-stable if and only if both $\Alg,\Alg'$ are rank-stable.
\end{theorem}

We prove this theorem in \Pref{prop:dir-prod-if} and \Pref{prop:dir-prod-only-if}.

\begin{lemma}\label{lem:idemp}
    Let $C\in\M[n](F)$ be a matrix such that $\rank(C^2-C)\le\lambda$. Then there exists a matrix $E\in\M[n](F)$ such that $E^2=E$ and $\rank(E-C)\le\lambda$.
\end{lemma}

\begin{proof}
    Take a basis $\set{v_1,\dots,v_r}$ of $\ker(C-\Id_n)$, and extend it to a basis $\set{v_1,\dots,v_k}$ of $\ker(C-\Id_n)\oplus\ker C$ by adding elements of $\ker C$ if necessary. Then, extend the resulting set to a basis $\set{v_1,\dots,v_n}$ of $F^n$. By the assumption, $\Rank(C^2-C)\le\lambda n$, so 
    \begin{equation} \label{eq:k n}
    k = \dim_F\ker C+\dim_F\ker(C-\Id_n)\ge\dim_F\ker(C^2-C)\ge(1-\lambda)n. 
   \end{equation}
    Put the vectors $v_1,\dots,v_n$ as the columns of a matrix $D\in\GL[n](F)$. We claim that
    \[
        D^{-1}CD=\begin{pmatrix}
            \Id_r\oplus 0_{k-r} & * \\ 0 & *
        \end{pmatrix}.
    \]
    Indeed, let $1\le i\le k$. If $i\le r$, then $D^{-1}CDe_i=D^{-1}Cv_i=D^{-1}v_i=e_i$; if $r+1\le i\le k$, then $D^{-1}CDe_i=D^{-1}Cv_i=0$. Therefore $D^{-1}CD$ has the above form.
    Take $E=D(\Id_r\oplus 0_{n-r})D^{-1}$. It is straightforward that $E^2=E$, and we have
    \[
        \Rank(E-C)=\Rank(D^{-1}ED-D^{-1}CD)\le n-k\le \lambda n,
    \]
    (the last inequality follows from \eqref{eq:k n}), so the lemma follows.
\end{proof}

\begin{lemma}\label{lem:idemp-block}
    Let $E\in\M[n](F)$ be a matrix of rank $r = \Rank E$ such that $E^2=E$. Then there exists an invertible matrix $D\in\GL[n](F)$ such that $DED^{-1}=\Id_r\oplus 0_{n-r}$. Furthermore, if $M\in\M[n](F)$ is a matrix satisfying
    \begin{align}
        \rank(ME-M) & \le \lambda \label{eq:idemp1}\\
        \rank(EM-M) & \le \lambda,\label{eq:idemp2}
    \end{align}
    then there exists a matrix $M^{\circ}\in\M[r](F)$ such that for $M'=D^{-1}(M^{\circ}\oplus 0_{n-r})D$ we have
    \[
        \rank(M-M') \le 2\lambda.
    \]
\end{lemma}

\begin{proof}
    The first part follows since every idempotent matrix is diagonalizable. We therefore prove only the second part. Let $M\in\M[n](F)$ be a matrix as above. Write the matrix $DMD^{-1}$ as a block matrix
    \[
        DMD^{-1}=\begin{pmatrix}
            M_{11} & M_{12} \\ M_{21} & M_{22}
        \end{pmatrix}
    \]
    for $M_{11}\in\M[r](F)$, $M_{12}\in\M[r\times (n-r)](F)$, $M_{21}\in\M[(n-r)\times r](F)$ and $M_{22}\in\M[n-r](F)$. Then
    \[
        D(ME-M)D^{-1}=\begin{pmatrix}
            0 & -M_{12} \\ 0 & -M_{22}
        \end{pmatrix}
    \]
    and
    \[
        D(EM-M)D^{-1}=\begin{pmatrix}
            0 & 0 \\ -M_{21} & -M_{22}
        \end{pmatrix}.
    \]
    By assumption, the rank of both of the above block triangular matrices is at most $\lambda n$. Let $M^{\circ}=M_{11}\in\M[r](F)$ and $M'=D^{-1}(M^{\circ}\oplus 0_{n-r})D$. Then
    \[
        D(M-M')D^{-1}=\begin{pmatrix}
            0 & M_{12} \\ M_{21} & M_{22}
        \end{pmatrix},
    \]
    and thus $\rank(M-M')\le 2\lambda$. The claim now follows.
\end{proof}

We turn to prove \Tref{thm:direct-product}.

\begin{proposition}\label{prop:dir-prod-if}
Let $\Alg,\Alg'$ be finitely presented algebras with finite-dimensional representations. If both $\Alg,\Alg'$ are rank-stable, then the direct product $\Alg\times \Alg'$ is rank-stable.
\end{proposition}

\begin{proof}
Fix presentations
\[
    \Alg=F\left<x_1,\dots,x_d\right>/\left<P_1,\dots,P_r\right>,\quad \Alg'=F\left<y_1,\dots,y_t\right>/\left<Q_1,\dots,Q_s\right>.
\]
Then the direct product $\Alg\times \Alg'$ has a presentation
\[
    \frac{F\left<x_1,\dots,x_d,y_1,\dots,y_t,e_1,e_2\right>}{\left<e_1P_j,\;e_2Q_j,\;e_1x_i-x_i,\;x_ie_1-x_i,\;e_2y_i-y_i,\;y_ie_2-y_i,\;e_1+e_2-1,\;e_1^2-e_1,\;e_2^2-e_2\right>}.
\] 
Let $L$ denote the maximal number of monomials in $P_j,Q_j$, and let $M$ denote the maximal length of a monomial appearing in any  $P_j$ or $Q_j$. Also, let $\rho\colon\Alg\to\M[p](F)$ and $\rho'\colon\Alg'\to\M[p'](F)$ be finite-dimensional representations of $\Alg$ and $\Alg'$.

We assume that $\Alg,\Alg'$ are rank-stable. Let $\eps>0$. Let $\delta_1$ and $\delta_2$ be the $\delta$'s from the rank-stability assumptions of $\Alg$ and $\Alg'$ with respect to $\eps$, respectively. We may assume that $\delta_1\le (\eps_1 d)^2$ and $\delta_2\le (\eps_2 t)^2$. Take $\delta>0$ small enough (during the proof, we will give upper bounds on $\delta$ in terms of the other parameters, which are set before $\delta$). Let $A_1,\dots,A_d,B_1,\dots,B_t,C_1,C_2\in\M[n](F)$ be matrices that $\delta$-approximate the relations from the above presentation of $\Alg\times \Alg'$.

First, note that $\rank(C_1^2-C_1)<\delta$ and $\rank(C_1+C_2-\Id_n)<\delta$. By \Lref{lem:idemp}, there exists an idempotent matrix $E_1\in\M[n](F)$ such that $\rank(C_1-E_1)<\delta$. Letting $E_2=\Id_n-E_1$, we have $E_2^2=E_2$, and
\[
    \rank(E_2-C_2)=\rank(\Id_n-E_1-C_2)\le\rank(C_1-E_1)+\rank(\Id_n-C_1-C_2)< 2\delta.
\]
Therefore the matrices $E_1,E_2$ approximate $C_1,C_2$. Denote $k=\Rank E_1$ and observe that there is an invertible matrix $D\in \GL[n](F)$ such that  $E_1=D^{-1}(\Id_k\oplus 0_{n-k})D$ and $E_2=D^{-1}(0_k\oplus \Id_{n-k})D$.

Next, for any $1\le i\le d$ we have $\rank(C_1A_i-A_i)<\delta$ and $\rank(A_iC_1-A_i)<\delta$. Hence:
\[
    \rank(E_1A_i-A_i)\le\rank(C_1A_i-A_i)+\rank(C_1-E_1)\le 2\delta
\]
and similarly $\rank(A_iE_1-A_i)\le 2\delta$. By \Lref{lem:idemp-block}, there exists a matrix $A_i^{\circ}\in\M[k](F)$ such that, for $A_i'=D^{-1}(A_i^{\circ}\oplus 0_{n-k})D$, we have $\rank(A_i-A_i')<4\delta$. The assumptions on $A_1,\dots,A_d$ also guarantee that $\rank (E_1P_j(A_1,\dots,A_d))<\delta$. Hence, by \Lref{lem:computing polyrank}, we have
\begin{align*}
    \rank(E_1P_j(A_1',\dots,A_d')) & \le \rank(E_1P_j(A_1,\dots,A_d)) \\
    & \quad + \rank(E_1P_j(A_1',\dots,A_d') -E_1P_j(A_1,\dots,A_d)) \\
    & \le \delta + \rank(P_j(A_1',\dots,A_d')-P_j(A_1,\dots,A_d)) \\
    & \le \delta+4LM\delta = (1+4LM)\delta.
\end{align*}
Note that:
\begin{align*}
    E_1P_j(A_1',\dots,A_d') &= D^{-1}(\Id_k\oplus 0_{n-k})(P_j(A_1^{\circ},\dots,A_d^{\circ})\oplus P_j(0_{n-k},\dots,0_{n-k}))D \\
    & = D^{-1}(P_j(A_1^{\circ},\dots,A_d^{\circ})\oplus 0_{n-k})D,
\end{align*}
hence
\[
    \Rank P_j(A_1^{\circ},\dots,A_d^{\circ}) = \Rank (E_1P_j(A_1',\dots,A_d')) \le (1+4LM)\delta n.
\]

For any $1\le i\le t$, we have 
\[
    \rank(E_2B_i-B_i)\le\rank(C_2B_i-B_i)+\rank(C_2-E_2)\le 3\delta
\]
and similarly $\rank(B_i E_2 - B_i)\le 3\delta$. Hence, for any $1\le i\le t$ there exists a matrix $B_i^{\circ}\in\M[n-k](F)$ such that, for $B_i'=D^{-1}(0_k\oplus B_i^{\circ})D$, we have $\rank(B_i-B_i')<6\delta$ and $\Rank Q_j(B_1^{\circ},\dots,B_t^{\circ})<(1+6LM)\delta n$.

We may assume that $\delta$ is small enough so that, if $A_1,\dots,A_d,B_1,\dots,B_t,C_1,C_2$ are not an exact solution for the relations of $\Alg\times\Alg'$, then $n\ge \max\set{d+\frac{p}{\eps},t+\frac{p'}{\eps}}$.

\begin{lemma}
    There exists a positive integer $k\le k_1\le k+\eps dn$ and matrices $A_1^{\bullet},\dots,A_d^{\bullet}\in\M[k_1](F)$ such that $P_j(A_1^{\bullet},\dots,A_d^{\bullet})=0_{k_1}$ for all $1\le j\le r$, and such that $\Rank(\widehat{A_i^{\circ}}-\widehat{A_i^{\bullet}}) \le (2+d)\eps n$ for all $1\le i\le d$.
\end{lemma}

\begin{proof}
    We divide into two cases.

    CASE I: $k\le \eps n$. Let $m=\lceil{\frac{k}{p}}\rceil$, and take $k_1=mp\le k+p$ and $A_i^{\bullet}=\Id_m\otimes\rho(x_i)$. Then $P_j(A_1^{\bullet},\dots,A_d^{\bullet})=0_{k_1}$ for all $1\le j\le r$, since $\rho$ is a representation of $\Alg$. So
    \[
        \Rank(\widehat{A_i^{\circ}}-\widehat{A_i^{\bullet}}) \le k+mp \le 2\eps n+p \le 3\eps n.
    \]

    CASE II: $k>\eps n$. Then
    \[
        \Rank P_j(A_1^{\circ},\dots,A_d^{\circ}) \le (1+4LM)\delta n \le \frac{(1+4LM)\delta}{\eps}k.
    \]
    For a small enough $\delta$ we may assume $\frac{(1+4LM)\delta}{\eps}<\delta_1$. Now we may use the rank-stability of $\Alg$; this proves that there exist $k''\ge 1$ and $A_1'',\dots,A_d''\in\M[k''](F)$ such that $P_j(A_1'',\dots,A_d'')=0_{k''}$ for all $1\le j\le r$, and such that $\Rank(\widehat{A_i^{\circ}}-\widehat{A_i''})\le\eps k$. Also, by \Lref{lem:small-approx}, we may assume that $(1-\eps d)k\le k''\le (1+\eps d)k$. If $k''\ge k$, the lemma is proved. Otherwise, take $m=\lceil \frac{k-k''}{p}\rceil$, and then the matrices $A_i^{\bullet}=A_i''\oplus (\Id_m\otimes\rho(x_i))$ of dimension $k_1=k''+mp\le k+p$ satisfy $P_j(A_1^{\bullet},\dots,A_d^{\bullet})=0_{k_1}$ for all $j$, and
    \begin{align*}
        \Rank(\widehat{A_i^{\circ}}-\widehat{A_i^{\bullet}}) &\le \Rank(\widehat{A_i^{\circ}}-\widehat{A_i''}) + \Rank(\widehat{A_i''}-\widehat{A_i^{\bullet}}) \\
        & \le \eps k + mp \le \eps n + k - k'' + p \le \eps n + \eps d k + \eps n \le (2+d)\eps n.
    \end{align*}
Thus proving the claim of the lemma.
\end{proof}

We now apply the above lemma to both $A_i^{\circ}$ and $B_i^{\circ}$: There exist $k\le k_1\le k+\eps d n$ and $n-k\le k_2\le n-k+\eps t n$, and matrices $A_1^{\bullet},\dots,A_d^{\bullet}\in\M[k_1](F)$ and $B_1^{\bullet},\dots,B_t^{\bullet}\in\M[k_2](F)$ such that $P_j(A_1^{\bullet},\dots,A_d^{\bullet})=0_{k_1}$ for all $j$ and $Q_j(B_1^{\bullet},\dots,B_t^{\bullet})=0_{k_2}$ for all $j$, and such that
\begin{align*}
    \Rank(\widehat{A_i^{\circ}}-\widehat{A_i^{\bullet}})& \le(2+d)\eps n \\
    \Rank(\widehat{B_i^{\circ}}-\widehat{B_i^{\bullet}})& \le(2+t)\eps n.
\end{align*}

Let $n'=k_1+k_2$, so that $n\le n'\le (1+\eps d+\eps t)n$. Let $D'\in\GL[n'](F)$ be the permutation matrix with rows
\[
    e_1,\dots,e_k,e_{k_1+1},\dots,e_{k_1+(n-k)},e_{k+1},\dots,e_{k_1},e_{k_1+(n-k)+1},\dots,e_{k_1+k_2}.
\]
This matrix satisfies the following property: for any $X\in\M[n'](F)$, the matrix $D'X(D')^{-1}$ is obtained from $X$ by moving the rows (resp.\ columns) $k+1,\dots,k_1$ to appear as rows $n+1,\dots,n+k_1-k$ (resp.\ columns). Also let $\widetilde{D}=D^{-1}\oplus\Id_{n'-n}$.

We define matrices $A_1'',\dots,A_d'',B_1'',\dots,B_t'',E_1'',E_2''\in\M[n'](F)$ by
\begin{align*}
    A_i'' &= \widetilde{D}D'(A_i^{\bullet}\oplus 0_{k_2})(\widetilde{D}D')^{-1},\\
    B_i'' &= \widetilde{D}D'(0_{k_1}\oplus B_i^{\bullet})(\widetilde{D}D')^{-1},\\
    E_1'' &= \widetilde{D}D'(\Id_{k_1}\oplus 0_{k_2})(\widetilde{D}D')^{-1},\\
    E_2'' &= \widetilde{D}D'(0_{k_1}\oplus \Id_{k_2})(\widetilde{D}D')^{-1}=\Id_{n'}-E_1'',
\end{align*}
and claim that they are the desired approximation to the original matrices. By their construction, it is immediate that $A_1'',\dots,A_d'',B_1'',\dots,B_t'',E_1'',E_2''$ are an exact solution to the relations of $\Alg\times\Alg'$. It remains to show that they well-approximate the matrices $A_1,\dots,A_d,B_1,\dots,B_t,C_1,C_2$.

First, we have
\[
    \widehat{\widetilde{D}^{-1}}\widehat{A_i'}\widehat{\widetilde{D}}=\widehat{DA_i'D^{-1}}=\widehat{A_i^{\circ}}
\]
and
\[
    \widehat{D'}\widehat{A_i^{\circ}} = \widehat{A_i^{\circ}}\widehat{D'} = \widehat{A_i^{\circ}},
\]
since $D'$ acts as the identity on the first $k\times k$ submatrix. Thus,
\begin{align*}
    \Rank(\widehat{A_i}-\widehat{A_i''}) & \le \Rank(\widehat{A_i}-\widehat{A_i'}) + \Rank(\widehat{A_i'}-\widehat{A_i''}) \\
    & = \Rank(A_i-A_i') + \Rank(\widehat{(\widetilde{D}D')^{-1}}\widehat{A_i'}\widehat{\widetilde{D}D'}-\widehat{(\widetilde{D}D')^{-1}}\widehat{A_i''}\widehat{\widetilde{D}D'})  \\
    & = \Rank(A_i-A_i') + \Rank(\widehat{A_i^{\circ}}-\widehat{A_i^{\bullet}})  \\
    & \le 4\delta n+ (2+d)\eps n \le (3+d)\eps n
\end{align*}
for a small enough $\delta$.

For $B_1,\dots,B_t$, note that
\[
    \widehat{\widetilde{D}^{-1}}\widehat{B_i'}\widehat{\widetilde{D}}=\widehat{DB_i'D^{-1}}=\widehat{0_k\oplus B_i^{\circ}}
\]
while
\[
    \widehat{(D')^{-1}}\widehat{(0_k\oplus B_i^{\circ})}\widehat{D'} = \widehat{(0_{k_1}\oplus B_i^{\circ})},
\]
so
\begin{align*}
    \Rank(\widehat{B_i}-\widehat{B_i''}) & \le \Rank(\widehat{B_i}-\widehat{B_i'}) + \Rank(\widehat{B_i'}-\widehat{B_i''}) \\
    & = \Rank(B_i-B_i') + \Rank(\widehat{(\widetilde{D}D')^{-1}}\widehat{B_i'}\widehat{\widetilde{D}D'}-\widehat{(\widetilde{D}D')^{-1}}\widehat{B_i''}\widehat{\widetilde{D}D'})  \\
    & = \Rank(B_i-B_i') + \Rank(\widehat{0_{k_1}\oplus B_i^{\circ}}-\widehat{0_{k_1}\oplus B_i^{\bullet}})  \\
    & \le 6\delta n+ (2+t)\eps n \le (3+t)\eps n
\end{align*}
for a small enough $\delta$. Finally,
\begin{align*}
    \Rank(\widehat{C_1}-\widehat{E_1''}) & \le \Rank(\widehat{C_1}-\widehat{E_1}) + \Rank(\widehat{E_1}-\widehat{E_1''}) \\
    & = \Rank(C_1-E_1) + \Rank(\widehat{(\widetilde{D}D')^{-1}}\widehat{E_1}\widehat{(\widetilde{D}D')^{-1}}-\widehat{(\widetilde{D}D')^{-1}}\widehat{E_1''}\widehat{(\widetilde{D}D')^{-1}}) \\
    & = \Rank(C_1-E_1) + \Rank(\widehat{\Id_k}-\widehat{\Id_{k_1}}) \\
    & \le \delta n + (k_1-k) \le \delta n + \eps d n \le (1+d)\eps n,
\end{align*}
and similarly $\Rank(\widehat{C_2}-\widehat{E_2''})\le (1+t)\eps n$ for small enough $\delta$. This proves that $\Alg\times\Alg'$ is rank-stable.
\end{proof}

\begin{proposition}\label{prop:dir-prod-only-if}
    Let $\Alg,\Alg'$ be finitely presented algebras with finite-dimensional representations. If $\Alg\times\Alg'$ is rank-stable, then both of $\Alg,\Alg'$ are rank-stable.
\end{proposition}

\begin{proof}

We use the same presentations of $\Alg$, $\Alg'$ and $\Alg\times\Alg'$ as in the proof of \Pref{prop:dir-prod-if}. Let $B_1,\dots,B_d\in\M[u](F)$ be an exact solution to $P_1,\dots,P_r$ and let $C_1,\dots,C_t\in\M[u'](F)$ be an exact solution to $Q_1,\dots,Q_s$, which exists since $\Alg,\Alg'$ have finite-dimensional representations. Let $\eps>0$ (we may readily assume that $\eps<1$), and let $\delta>0$ be the $\delta$ from the rank-stability of $\Alg\times\Alg'$ with respect to $\eps'=\frac{\eps}{2(d+t+6)}$. We assume that $\delta<(\eps'(d+t+2))^2$. Take $\delta'<\min\big\{\frac{\delta}{2},\frac{\eps'}{u},\frac{\eps'}{u'}\big\}$.

Let $A_1,\dots,A_d\in\M[n](F)$ be matrices such that
\[
    \rank P_j(A_1,\dots,A_d)<\delta'.
\]
Since $\delta'<\frac{\eps'}{u'}$, if $n\le u'/\eps'$ then $A_1,\dots,A_d$ is an exact solution to $P_1,\dots,P_r$, and we are done. We therefore assume that $\eps'n>u'$. Similarly, we may assume $\eps'n>u$. Define matrices $A_1^{\circ},\dots,A_d^{\circ},C_1^{\circ},\dots,C_t^{\circ}\in\M[n+u'](F)$ by $A_i^{\circ}=A_i\oplus 0_{u'}$ and $C_i^{\circ}=0_n\oplus C_i$. Define in addition $E_1=\Id_n\oplus 0_{u'}\in\M[n+u'](F)$ and $E_2=\Id_{n+u'}-E_1=0_n\oplus \Id_{u'}$. It follows that
\begin{align*}
    \Rank (E_1P_j(A_1^{\circ},\dots,A_d^{\circ}))&=\Rank(P_j(A_1,\dots,A_d))\le\delta' n<\delta(n+u'),\\
    E_2Q_j(C_1^{\circ},\dots,C_d^{\circ})&=0_n\oplus Q_j(C_1,\dots,C_t)=0_{n+u'}\\
    E_1A_i^{\circ}&=A_i^{\circ}E_1=A_i^{\circ}\\
    E_2C_i^{\circ}&=C_i^{\circ}E_2=C_i^{\circ}\\
    E_1^2&=E_1\\
    E_2^2&=E_2.
\end{align*}
Using the rank-stability of $\Alg\times \Alg'$, there exist an integer $n'\ge 1$ and matrices $A_1',\dots,A_d',\allowbreak C_1',\dots,C_t',E_1',E_2'\in\M[n'](F)$ that satisfy the relations of $\Alg\times\Alg'$, and $\eps'$-approximate $A_1^{\circ},\dots,A_d^{\circ},C_1^{\circ},\dots,C_t^{\circ},E_1,E_2$. Also, since $\delta<(\eps'(d+t+2))^2$, we may assume 
\[
\left|n'-(n+u')\right|\le\eps'(d+t+2)(n+u').
\]
In particular,
\[
    n'\le (1+\eps'(d+t+2))(n+u') \le (1+\eps'(d+t+2))(1+\eps') n.
\]

Note that $(E_1')^2=E_1'$ and $E_1'A_i'=A_i'E_1'=A_i'$ for all $1\le i\le d$. It follows that there exist an invertible matrix $D\in\GL[n'](F)$ and $0\le k\le n'$ such that $DE_1'D^{-1}=\Id_k\oplus 0_{n'-k}$. Also, for any $1\le i\le d$ we have $DA_i'D^{-1}=A_i''\oplus 0_{n'-k}$ for some $A_i''\in\M[k](F)$. Note that $A_1'',\dots,A_d''$ are an exact solution to $P_1,\dots,P_r$, because
\begin{align*}
    0_{n'} = E_1'P_j(A_1',\dots,A_d') &= D^{-1}((\Id_k\oplus 0_{n'-k})(P_j(A_1'',\dots,A_d'')\oplus P_j(0,\dots,0)))D \\
    & =D^{-1}(P_j(A_1'',\dots,A_d'')\oplus 0_{n'-k})D
\end{align*}
for all $1\le j\le r$. Also, since $E_1'$ is an $\eps'$-approximation of $E_1=\Id_n\oplus 0_{u'}$,
\[
    \left|n-k\right|=\left|\Rank(\widehat{E_1})-\Rank(\widehat{E_1'})\right|\le \Rank(\widehat{E_1}-\widehat{E_1'})\le\eps'n,
\]
so we must have $(1-\eps')n\le k\le(1+\eps')n$.

Let $n''$ be the unique integer satisfying $\max\set{n,n'}\le n''<\max\set{n,n'}+u$ and $n''\equiv k\pmod{u}$. Write $n''=k+qu$, and let $D'=D\oplus\Id_{n''-n'}\in\GL[n''](F)$. For every $1\le i\le d$, define a matrix $M_i\in\M[n''](F)$ by $M_i=(D')^{-1}(A_i''\oplus (\Id_q\otimes B_i))D'$. Since $A_1'',\dots,A_d''$ and $B_1,\dots,B_d$ are exact solutions of $P_1,\dots,P_r$, so are $M_1,\dots,M_d$. Also,
\begin{align*}
    \Rank(\widehat{M_i}-\widehat{A_i}) & \le \Rank(\widehat{M_i}-\widehat{A_i'}) + \Rank(\widehat{A_i'}-\widehat{A_i}) \\
    & = \Rank(\widehat{D'}\widehat{M_i}\widehat{(D')^{-1}} - \widehat{D'}\widehat{A_i'}\widehat{(D')^{-1}}) + \Rank(\widehat{A_i'}-\widehat{A_i}) \\
    & = \Rank((A_i''\oplus (\Id_q\otimes B_i)) - (A_i''\oplus 0_{qu})) + \Rank(\widehat{A_i'}-\widehat{A_i}) \\
    & \le qu + \eps'n = (n''-k) + \eps'n \le \max\set{n,n'}+u-k+\eps'n \\
    & \le (1+\eps'(d+t+2))(1+\eps') n + \eps'n - (1-\eps')n + \eps'n \\
    & = ((d+t+6)\eps'+(d+t+2)(\eps')^2)n \le \eps n.
\end{align*}
This shows that $\Alg$ is rank-stable, as required.
\end{proof}

\subsection{Matrix algebras}

\begin{theorem}\label{thm:mats}
    Let $\Alg$ be a finitely presented $F$-algebra that has a finite-dimensional representation, and let $m\ge 1$. Then $\Alg$ is rank-stable if and only if $\M[m](\Alg)$ is rank-stable.
\end{theorem}

For the rest of the section, we denote by $e_{ij}$ the standard matrix units of size $m\times m$. Before proving the theorem, we give the following useful lemma:
\begin{lemma}\label{lem:new8.2}
    Let $\set{E_{ij}}_{i,j=1}^m\sub\M[n'](F)$ be a set of matrices satisfying
    \[
        E_{ij}E_{kl}=\delta_{jk}E_{il} \ \ \textrm{and}\ \ \sum_{i=1}^m E_{ii}=\Id_{n'}.
    \]
    Then $n'=qm$ for some integer $q\ge 1$, and there exists an invertible matrix $D\in\GL[n'](F)$ such that $DE_{ij}D^{-1}=e_{ij}\otimes \Id_q$. Furthermore, let $n,\lambda\ge 1$ be integers, and set $N=\left|n'-n\right|$. If $A\in\M[n](F)$ is a matrix satisfying
    \begin{equation}\label{eq:A-E_ij}
        \Rank(\widehat{A}\widehat{E_{ij}}-\widehat{E_{ij}}\widehat{A})\le\lambda
    \end{equation}
    for all $1\le i,j\le d$, then there exists a matrix $A^{\circ}\in\M[q](F)$ such that for the matrix $C=D^{-1}(\Id_m\otimes A^{\circ})D\in\M[n'](F)$ we have
    \[
        \Rank(\widehat{A}-\widehat{C})\le m^2(\lambda+4N).
    \]
\end{lemma}

\begin{proof}
    The first part is classical, and we briefly recall its proof. The matrices $\set{E_{ij}}$ form a set of $m\times m$ matrix units, hence induce a unital embedding $\M[m](F)\hookrightarrow\M[n'](F)$. It follows that $n'=qm$ for some integer $q\ge 1$. In addition, by the Skolem--Noether theorem, any such embedding must be conjugate to the embedding $X\mapsto X\otimes \Id_q$ via some invertible matrix $D\in\GL[n'](F)$, and therefore $DE_{ij}D^{-1}=e_{ij}\otimes \Id_q$ for all $1\le i,j\le m$.

    For the second part, let $A\in\M[n](F)$ be a matrix satisfying \eqref{eq:A-E_ij} for all $1\le i,j\le m$. Let $A'\in\M[n'](F)$  be the matrix obtained from $A$ by adding zero rows and columns if $n\le n'$, or removing the last $n-n'$ rows and columns if $n>n'$. It follows that
    \begin{equation}\label{eq:A-A'}
        \Rank(\widehat{A}-\widehat{A'})\le 2N.
    \end{equation}

    We work with the matrix $A'$. Substituting in \eqref{eq:A-E_ij}, we have
    \begin{align*}
        \Rank(A'E_{ij}-E_{ij}A') & = \Rank(\widehat{A'}\widehat{E_{ij}}-\widehat{E_{ij}}\widehat{A'}) \\
        & \le \Rank(\widehat{A}\widehat{E_{ij}}-\widehat{E_{ij}}\widehat{A}) + 2 \Rank(\widehat{A}-\widehat{A'}) \\
        & \le \lambda + 4N
    \end{align*}
    for all $1\le i,j\le m$.

    Write $DA'D^{-1}$ as an $m\times m$ block matrix
    \[
        DA'D^{-1}=\begin{pmatrix}
            A'(1,1) & \cdots & A'(1,m) \\ \vdots & \ddots & \vdots \\ A'(m,1) & \cdots & A'(m,m)
        \end{pmatrix}
    \]
    where each $A'(i,j)\in\M[q](F)$. Then the matrix $D(E_{ij}A'-A'E_{ij})D^{-1}$, as an $m\times m$ block matrix, contains only blocks in its $i$-th row and $j$-th column of blocks. For $i\ne j$, the block that appears in the $(i,i)$ position of $D(E_{ij}A'-A'E_{ij})D^{-1}$ is precisely $A'(j,i)$, whereas the block that appears in the $(i,j)$ position is $A'(j,j)-A'(i,i)$. It follows that for all $i\ne j$ we have
    \[
        \Rank(A'(j,i)) \le \lambda + 4N
    \]
    and
    \[
        \Rank(A'(j,j)-A'(i,i)) \le \lambda + 4N.
    \]

    Write $A^{\circ}=A'(1,1)$ and $C=D^{-1}(\Id_m\otimes A^{\circ})D$. The matrix $DA'D^{-1}-DCD^{-1}$ has the block matrix form
        \[
            \begin{pmatrix}
                A'(1,1) - A'(1,1) & A'(1,2) & \cdots & A'(m,m) \\ 
                A'(2,1) & A'(2,2) - A'(1,1) &  & \vdots \\
                \vdots &  & \ddots & A'(m-1,m) \\
                A'(m,1) & \cdots & A'(m,m-1) & A'(m,m) - A'(1,1)
            \end{pmatrix}
        \]
        In other words, the matrix $DA'D^{-1}-DCD^{-1}$, as a block matrix, contains the blocks $A'(i,j)$ off-diagonal, and $A'(i,i)-A'(1,1)$ along the main diagonal. Therefore
        \[
            \Rank(A'-C)=\Rank(DA'D^{-1}-DCD^{-1})\le (m^2-1)(\lambda+4N).
        \]
        We may now conclude that
        \begin{align*}
            \Rank(\widehat{A}-\widehat{C}) & \le \Rank(\widehat{A}-\widehat{A'}) + \Rank(\widehat{A'}-\widehat{C}) \\
            & = \Rank(\widehat{A}-\widehat{A'}) + \Rank(A'-C) \\
            & \le 2N + (m^2-1)(\lambda + 4N) \le m^2(\lambda+4N)
        \end{align*}
        as required.
\end{proof}

We now prove \Tref{thm:mats} in two parts.

\begin{proposition}
    If $\Alg$ is a finitely presented rank-stable $F$-algebra that admits a finite-dimensional representation, then $\M[m](\Alg)$ is rank-stable for every $m\ge 1$.
\end{proposition}

\begin{proof}
    Let $\Alg$ be a finitely presented $F$-algebra. Fix a presentation
    \[
        \Alg=F\left<x_1,\dots,x_d\right>/\left<P_1,\dots,P_r\right>.
    \]
    Let $L$ denote the maximal number of monomials in $P_1,\dots,P_r$, and $M$ the maximal degree of $P_1,\dots,P_r$. Also, fix a solution $M_1,\dots,M_d\in\M[s](F)$ to $P_1=\cdots=P_r=0_s$, which exists since $\Alg$ has a finite-dimensional representation. Recall that $\M[m](F)$ admits a presentation
    \[
        \M[m](F)=F\left<e_{ij}\suchthat 1\le i,j\le m\right>/\left<e_{ij}e_{kl}-\delta_{jk}e_{il},\;{\textstyle \sum_i} e_{ii}-1\right>.
    \]
    Note that there is a slight abuse of notations here\mdash{}we use $e_{ij}$ to denote both abstract generators of $\M[m](F)$, and the explicit matrix units in $\M[m](F)$. Then the algebra $\M[m](\Alg)=\Alg\otimes_F \M[m](F)$ has the presentation
    \[
        \M[m](\Alg) = \frac{F\left<x_k,e_{ij}|1\leq k\leq d,1\leq i,j\leq m\right>}{\left<P_1,\dots,P_r,\; e_{ij}e_{kl} - \delta_{jk}e_{il},\; {\textstyle \sum_i} e_{ii} - 1,\;e_{ij}x_k - x_k e_{ij}\right>}.
    \]

    Before giving the proof, we sketch the idea behind it. Suppose we have matrices $A_k,B_{ij}$ that are approximate solutions of the relations of $\M[m](\Alg)$ in the above presentation. The matrices $B_{ij}$ are approximate solutions of the relations of $\M[m](F)$, hence are close to a solution to these relations (by \Tref{thm:fin-dim}), which is\mdash{}up to conjugation\mdash{}a system of block matrix units $\{E_{ij}\}_{i,j=1}^{m}$. The matrices $A_k$ almost commute with every $E_{ij}$, hence are close to being block scalar matrices. The obtained blocks are an approximate solution to the relations of $\Alg$, hence are close to a solution to these relations, by the stability assumption on $\Alg$. Tensoring these blocks with $E_{ij}$, we get the desired approximation.

    Recall that $\Alg$ is rank-stable. We want to prove that $\M[m](\Alg)$ is rank-stable. Let $\eps>0$. We make the following choices of constants, in this order: we fix $\eps_1>0$, which depends only on $\eps$; we then choose $\delta_1>0$ to be the $\delta$ from the rank-stability of $\Alg$ with respect to $\eps_1$, and which will depend on $\eps,\eps_1$; next, we choose $\eps_2>0$ small enough that depends on $\eps,\eps_1,\delta_1$, and then $\delta_2>0$ to be the $\delta$ from the rank-stability of $\M[m](F)$ with respect to $\eps_2$. Finally, we choose $\eta>0$, sufficiently small depending on all previous constants. Throughout the proof we will add constraints on these constants, which will be compatible with the above dependencies.

    Let $A_1,\dots,A_d,B_{11},B_{12},\dots,B_{mm}\in\M[n](F)$ be matrices such that
    \begin{align}
        \rank(P_j(A_1,\dots,A_d)) & < \eta \label{eq:Ak-Ps} \\
        \rank(B_{ij}B_{kl}-\delta_{jk}B_{il}) & < \eta \label{eq:Bij-mult}\\
        \rank({\textstyle \sum_i} B_{ii}-\Id_n) & < \eta \label{eq:Bii-sum}\\
        \rank(B_{ij}A_k-A_kB_{ij}) & < \eta. \label{eq:Ak-Bij}
    \end{align}
    We want to $\eps$-approximate these matrices by a solution to the above presentation of $\M[m](\Alg)$.

    First, note that by \eqref{eq:Bij-mult} and \eqref{eq:Bii-sum}, the matrices $\vec{B}=(B_{ij})$ are an approximate solution to the relations of $\M[m](F)$. Assuming $\eta<\delta_2$, there exist $n'\ge 1$ and matrices $E_{ij}\in\M[n'](F)$ such that
    \begin{align*}
        E_{ij}E_{kl} &= \delta_{jk}E_{il},\\
        \sum_{i=1}^m E_{ii} &= \Id_{n'},\\
        \Rank(\widehat{B_{ij}}-\widehat{E_{ij}}) &< \eps_2 n.
    \end{align*}
    Also, by assuming that $\eta<\delta_2<(\eps_2m^2)^2$, \Lref{lem:small-approx} shows that we may choose $n'$ such that
    \[
        (1-\eps_2m^2)n\le n'\le (1+\eps_2m^2)n,
    \]
    and so
    \begin{equation}\label{eq:nn'-mats}
        \left|n'-n\right|\le \eps_2m^2 n.
    \end{equation}

    By \eqref{eq:Ak-Bij}, it follows that for every $1\le k\le d$ and $1\le i,j\le m$,
    \begin{align*}
        \Rank(\widehat{E_{ij}}\widehat{A_k}-\widehat{A_k}\widehat{E_{ij}}) & \le \Rank(\widehat{B_{ij}}\widehat{A_k}-\widehat{A_k}\widehat{B_{ij}}) + 2\Rank(\widehat{B_{ij}}-\widehat{E_{ij}}) \\
        & =\Rank(B_{ij}A_k-A_kB_{ij}) + 2\Rank(\widehat{B_{ij}}-\widehat{E_{ij}}) \\
        & \le \eta n + 2\eps_2 n = (\eta+2\eps_2)n.
    \end{align*}
    We can now apply \Lref{lem:new8.2}. Thus, we can write $n'=qm$ for some integer $q\ge 1$; there exists an invertible matrix $D\in\GL[n'](F)$ such that $DE_{ij}D^{-1}=e_{ij}\otimes \Id_q$ for all $1\le i,j\le m$; and there exist matrices $A_1^{\circ},\dots,A_k^{\circ}\in\M[q](F)$ such that, writing $C_k=D^{-1}(\Id_m\otimes A_k^{\circ})D$, we have
    \begin{equation}\label{eq:Ak-Ck}
        \Rank(\widehat{A_k}-\widehat{C_k}) \le m^2((\eta+2\eps_2)n + 4\left|n'-n\right|)\le (\eta+6\eps_2)m^4n
    \end{equation}
    for all $1\le k\le d$.

    Next, we use \Lref{lem:computing polyrank}; combining \eqref{eq:Ak-Ps} and \eqref{eq:Ak-Ck}, we have
    \begin{align*}
        \Rank(P_j(C_1,\dots,C_d)) &\le \Rank(P_j(A_1,\dots,A_d)) \\
        & \quad + LM(\eta+6\eps_2)m^4n + \left|n'-n\right| \\
        & \le (\eta + LM(\eta+6\eps_2)m^4 + \eps_2 m^2) n \\
        & \le \frac{\eta + LM(\eta+6\eps_2)m^4 + \eps_2 m^2}{1-\eps_2m^2} n'.
    \end{align*}
    By choosing $\eta,\eps_2$ to be small enough, we may assume that
    \[
        \Rank(P_j(C_1,\dots,C_d))\le\delta_1 n'.
    \]
    By the definition of $C_1,\dots,C_d$, we have
    \[
        P_j(C_1,\dots,C_d)=D^{-1}(\Id_m\otimes P_j(A_1^{\circ},\dots,A_d^{\circ}))D,
    \]
    and so
    \[
        \Rank(P_j(A_1^{\circ},\dots,A_d^{\circ})) \le \delta_1 q.
    \]
    We now apply the rank-stability of $\Alg$ to find $q'\ge 1$ and matrices $A_1',\dots,A_k'\in\M[q'](F)$ such that
    \begin{equation}\label{eq:Ak'-Ps}
        P_j(A_1',\dots,A_d')=0_{q'}
    \end{equation}
    for all $1\le j\le r$, and such that
    \begin{equation}\label{eq:circ}
        \Rank(\widehat{A_k^{\circ}}-\widehat{A_k'})\le \eps_1 q
    \end{equation}
    for all $1\le k\le d$. By \Lref{lem:small-approx}, assuming $\delta_1<(\eps_1 d)^2$, we may also assume that
    \[
        (1-\eps_1 d)q\le q'\le (1+\eps_1 d)q.
    \]
    We note that we can assume that $q'\ge q$. Otherwise, let $t$ be the smallest integer such that $q'+ts\ge q$. Replacing each $A_k'$ by $A_k''\coloneqq A_k'\oplus (\Id_t\otimes M_k)\in\M[q'+ts](F)$, we still have
    \[
        P_j(A_1'',\dots,A_d'')=0_{q'+ts}
    \]
    for all $1\le j\le r$, and
    \[
        \Rank(\widehat{A_k^{\circ}}-\widehat{A_k''}) \le \Rank(\widehat{A_k^{\circ}}-\widehat{A_k'}) + ts \le \eps_1 q + |q'-q| + s \le \eps_1(d+1)q + s.
    \]
    
    Fixing $\eps,\eps_1,\delta_1,\eps_2,\delta_2$ and letting $\eta\rightarrow 0$, we see that either $A_1,\dots,A_d,B_{11},\dots,B_{mm}$ form an exact solution to the prescribed relations, or else $n\rightarrow \infty$, and so $n'\rightarrow \infty$ by \eqref{eq:nn'-mats} and so $q\rightarrow \infty$. Therefore, we can take $\eta>0$ to be sufficiently small, forcing $q\geq s/\eps_1$, ensuring that $\Rank(\widehat{A_k^{\circ}}-\widehat{A_k''})\leq \eps_1 (d+2) q$.

    Assume therefore $q'\ge q$, and work with $A_1',\dots,A_k'$. We define matrices
    \[
        \widetilde{A}_1,\dots,\widetilde{A}_k, \widetilde{B}_{11},\widetilde{B}_{12},\dots, \widetilde{B}_{mm}\in\M[mq'](F)
    \]
    as follows. Let $D'\in\M[mq'](F)$ denote the permutation matrix with rows: 
    \begin{multline*}
    e_1,\dots,e_q,e_{q'+1}, \dots ,e_{q'+q},e_{2q'+1},\dots,e_{2q'+q}, \dots, e_{(m-1)q'+q},
 \\ e_{q+1},\dots,e_{q'},e_{q'+q+1},\dots,e_{2q'},\dots, e_{mq'}.
    \end{multline*}
    In other words, for any matrix $X\in\M[mq'](F)$, the matrix $D'X(D')^{-1}$ is obtained from $X$ by moving the last $q'-q$ rows and columns of any of the $m$ blocks of size $q'\times q'$ of $X$ to the end of the matrix. Also, let $\widetilde{D}=D^{-1}\oplus \Id_{mq'-mq}$. For any $1\le k\le d$ and $1\le i,j\le m$ define \[ \widetilde{A}_k=\widetilde{D}D'(\Id_m\otimes A_k')(\widetilde{D}D')^{-1}\ \ \text{and}\ \ \widetilde{B}_{ij}=\widetilde{D}D'(e_{ij}\otimes \Id_{q'})(\widetilde{D}D')^{-1}. \] We claim that $\widetilde{A}_1,\dots,\widetilde{A}_d,\widetilde{B}_{11},\widetilde{B}_{12},\dots,\widetilde{B}_{mm}$ form a solution to the relations of $\M[m](\Alg)$, and $\eps$-approximate the matrices $A_1,\dots,A_d,B_{11},B_{12},\dots,B_{mm}$.

    Let us first prove that $\widetilde{A}_k,\widetilde{B}_{ij}$ are a solution to the relations of $\M[m](\Alg)$. The relations $P_j(\widetilde{A}_1,\dots,\widetilde{A}_d)=0_{mq'}$ follow from \eqref{eq:Ak'-Ps}; the relations $\widetilde{B}_{ij}\widetilde{B}_{kl}=\delta_{jk}\widetilde{B}_{il}$ and $\sum_{i=1}^m \widetilde{B}_{ii}=\Id_{mq'}$ follow from the definition of $\widetilde{B}_{ij}$; and it follows from their definition that $\widetilde{A}_k$ and $\widetilde{B}_{ij}$ commute. Therefore these matrices form a solution to the relations of $\M[m](\Alg)$.

    We now prove that this solution is $\eps$-close to the original approximate solution we started with. First, for every $i,j$, we can write:
    \[
    D'(e_{ij}\otimes\Id_{q'})(D')^{-1}=(e_{ij}\otimes\Id_q)\oplus (e_{ij}\otimes\Id_{q'-q}),
    \]
    so
    \[
        \widetilde{B}_{ij}=(D^{-1}(e_{ij}\otimes \Id_q)D)\oplus (e_{ij}\otimes\Id_{q'-q}) = E_{ij}\oplus(e_{ij}\otimes\Id_{q'-q}).
    \]
    Hence
    \begin{align*}
        \Rank(\widehat{B_{ij}}-\widehat{\widetilde{B}_{ij}}) & \le \Rank(\widehat{B_{ij}}-\widehat{E_{ij}}) + \Rank(\widehat{E_{ij}}-\widehat{\widetilde{B}_{ij}})\\
        & = \Rank(\widehat{B_{ij}}-\widehat{E_{ij}}) + (q'-q) \le \eps_2 n + \eps_1 dq \\
        & = \eps_2 n+\eps_1 d\frac{n'}{m} \le \eps_2 n + \eps_1 d(1+\eps_2 m^2)\frac{n}{m} \le  \eps n
    \end{align*}
    as long as we take $\eps_1,\eps_2$ small enough compared to $\eps$.
    
    Finally, fix $1\le k\le d$. The first $mq\times mq$ block of $D'(\Id_m\otimes A_k')(D')^{-1}$ is a block-scalar matrix, where the repeating block is the first $q\times q$ rows and columns of $A_k'$, denoted $[A_k']_q$. In other words, the first $mq\times mq$ block of $D'(\Id_m\otimes A_k')(D')^{-1}$ is $\Id_m \otimes [A_k']_q$. In particular,
    \begin{equation} \label{eq:Qatum}
    \Rank\left( \widehat{D'}\widehat{(\Id_m\otimes A_k')}\widehat{(D')^{-1}} - \widehat{\Id_n\otimes [A_k']_q} \right) \leq 2(mq' - mq).
    \end{equation}
    
    Hence:
    \begin{align*}
        \Rank(\widehat{C_k}-\widehat{\widetilde{A}_k}) &= \Rank(\widehat{D^{-1}}\widehat{(\Id_m\otimes A_k^{\circ})}\widehat{D}-\widehat{\widetilde{A}_k})\\
        &= \Rank(\widehat{\widetilde{D}}\widehat{(\Id_m\otimes A_k^{\circ})}\widehat{\widetilde{D}^{-1}}-\widehat{\widetilde{A}_k})\\
        &= \Rank\left(\widehat{\widetilde{D}}\widehat{(\Id_m\otimes A_k^{\circ})}\widehat{\widetilde{D}^{-1}}-\widehat{\widetilde{D}}\widehat{D'}\widehat{(\Id_m\otimes A_k')}\widehat{\widetilde{(D')^{-1}}}\widehat{\widetilde{D}^{-1}}\right)\\
        & = \Rank(\widehat{\Id_m\otimes A_k^{\circ}} - \widehat{D'}\widehat{(\Id_m\otimes A_k')}\widehat{(D')^{-1}}) \\
        &\le \Rank(\Id_m\otimes A_k^{\circ}-\Id_m\otimes [A_k']_q) + 2m(q'-q) \quad\qquad\text{by \eqref{eq:Qatum}}\\
        & = \Rank(\Id_m\otimes(A_k^{\circ}-[A_k']_q)) + 2m(q'-q) \\
        &= m\Rank(A_k^{\circ}-[A_k']_q) + 2m(q'-q) \\
        & \le m\Rank(\widehat{A_k^{\circ}}-\widehat{A_k'})+4m(q'-q) \\
        & \le m\eps_1 q+4m\eps_1 dq = (1+4d)\eps_1 n'.
    \end{align*}
    Summing up,
    \begin{align*}
        \Rank(\widehat{A_k}-\widehat{\widetilde{A}_k}) & \le \Rank(\widehat{A_k}-\widehat{C_k}) + \Rank(\widehat{C_k} - \widehat{\widetilde{A}_k})\\
        & \le (\eta+6\eps_2)m^4n + (1+4d)\eps_1 n' \le \eps n
    \end{align*}
    for a small enough choice of $\eps_1,\eps_2,\eta$ which only depends on $\eps$. This concludes the proof that if $\Alg$ is rank-stable, then so is $\M[m](\Alg)$.
\end{proof}

\begin{proposition}
    If $\Alg$ is finitely presented and $\M[m](\Alg)$ is rank-stable for some $m\ge 1$ then $\Alg$ is rank-stable.
\end{proposition}

\begin{proof}
    Suppose $\M[m](\Alg)$ is rank-stable. We want to prove that $\Alg$ is rank-stable. We use again the presentations
    \[
        \Alg=F\left<x_1,\dots,x_d\right>/\left<P_1,\dots,P_r\right>
    \]
    and
    \[
        \M[m](\Alg) = \frac{F\left<x_k,e_{ij}|1\leq k\leq d,1\leq i,j\leq m\right>}{\left<P_1,\dots,P_r,\; e_{ij}e_{kl} - \delta_{jk}e_{il},\; {\textstyle \sum_i} e_{ii} - 1,\;e_{ij}x_k - x_k e_{ij}\right>}.
    \]
    Let $\eps>0$, and let $\delta>0$ be the $\delta$ from the rank-stability of $\M[m](\Alg)$ with respect to $\eps'$, which will be taken to be small enough depending on $\eps,d,m$. We may assume that $\delta<(\eps'(d+m^2))^2$. Let $A_1,\dots,A_d\in\M[n](F)$ be matrices such that
    \[
        \rank(P_j(A_1,\dots,A_d))<\delta.
    \]

    Write $A_k'=\Id_m\otimes A_k$ and $B_{ij}=e_{ij}\otimes \Id_n$. It holds that
    \[
        \rank(P_j(A_1',\dots,A_d'))<\delta,
    \]
    and that the matrices $A_k'$ and $B_{ij}$ satisfy the other relations in the above presentation of $\M[m](\Alg)$. By the rank-stability of $\M[m](\Alg)$, there exist $n'\ge 1$ and $C_1,\dots,C_d,\allowbreak E_{11},E_{12},\dots,E_{mm}\in\M[n'](F)$ that satisfy the relations of $\M[m](\Alg)$, and such that
    \begin{equation}\label{eq:mats-Ak'-Ck}
        \Rank(\widehat{A_k'}-\widehat{C_k})\le\eps' mn
    \end{equation}
    and
    \begin{equation}\label{eq:mats-Bij-Eij}
        \Rank(\widehat{B_{ij}}-\widehat{E_{ij}})\le\eps' mn.
    \end{equation}
    By \Lref{lem:small-approx}, we may also assume that
    \begin{equation}\label{eq:mats-n'mn}
        (1-\eps'(d+m^2))mn\le n'\le (1+\eps'(d+m^2))mn.
    \end{equation}
    
    We apply the first part of \Lref{lem:new8.2}. Therefore $n'=qm$ for some integer $q\ge 1$, and there is an invertible matrix $D\in\GL[n'](F)$ such that $DE_{ij}D^{-1}=e_{ij}\otimes \Id_q$ for all $1\le i,j\le m$. Since $C_kE_{ij}=E_{ij}C_k$ for all $1\le i,j\le m$, there exists $C_k^{\circ}\in\M[q](F)$ such that $DC_kD^{-1}=\Id_m\otimes C_k^{\circ}$. Note that \eqref{eq:mats-n'mn} implies
    \[
        (1-\eps'(d+m^2))n\le q\le (1+\eps'(d+m^2))n,
    \]
    and since $C_1,\dots,C_d$ satisfy $P_1,\dots,P_r$ it follows that $C_1^{\circ},\dots,C_d^{\circ}$ satisfy these relations as well.

    Our next goal is to prove that $D$ is close to a block-diagonal matrix. Indeed, write $B_{ij}'=e_{ij}\otimes\Id_q\in\M[n'](F)$, so $DE_{ij}D^{-1}=B_{ij}'$. Note that for each $1\le i\le m$, the matrix $\widehat{B_{ii}}-\widehat{B_{ii}'}$ 
    is a diagonal matrix with at most $(2i-1)\left|n-q\right|$ non-zero entries. Hence
    \[
        \Rank\!\left(\widehat{B_{ii}}-\widehat{B_{ii}'}\right)\le (2i-1)\left|n-q\right| \le (2m-1)\left|n-q\right|\le 2\eps'(d+m^2)mn.
    \]
    It follows that for every $1\le i\le m$ we have
    \begin{align*}
        \Rank(DB_{ii}'-B_{ii}'D) & = \Rank(\widehat{D}\widehat{B_{ii}'}-\widehat{B_{ii}'}\widehat{D}) \\
        & \le \Rank(\widehat{D}\widehat{B_{ii}}-\widehat{B_{ii}'}\widehat{D}) + \Rank(\widehat{B_{ii}}-\widehat{B_{ii}'})\\
        & \le \Rank(\widehat{D}\widehat{E_{ii}}-\widehat{B_{ii}'}\widehat{D}) + \Rank(\widehat{B_{ii}}-\widehat{B_{ii}'}) + \Rank(\widehat{B_{ii}}-\widehat{E_{ii}})\\
        & = \Rank(DE_{ii}-B_{ii}'D) + \Rank(\widehat{B_{ii}}-\widehat{B_{ii}'}) + \Rank(\widehat{B_{ii}}-\widehat{E_{ii}}) \\
        & = \Rank(\widehat{B_{ii}}-\widehat{B_{ii}'}) + \Rank(\widehat{B_{ii}}-\widehat{E_{ii}}) \\
        & \le 2\eps'(d+m^2)mn+\eps'mn \le 3\eps'(d+m^2)mn.
    \end{align*}
    In a similar manner to the proof of the second part of \Lref{lem:new8.2}, it now follows that there exist matrices $D_1,\dots,D_m\in\M[q](F)$ such that, for $D'=D_1\oplus\cdots\oplus D_m\in\M[n'](F)$, we have
    \[
        \Rank(D-D') \le m^2(3\eps'(d+m^2)mn+4\left|n-q\right|) \le 7\eps'(d+m^2)m^3n.
    \]
    Since $D$ is invertible, it follows that $\sum_{i=1}^m\Rank(D_i)=\Rank(D')\ge n'-7\eps'(d+m^2)m^3n$. Hence, replacing each $D_i$ by an invertible matrix $D_i'\in\GL[q](F)$ such that $\Rank(D_i'-D_i)=q-\Rank(D_i)$, we may assume that $D_1,\dots,D_m$ are invertible, and $\Rank(D-D')\le 14\eps'(d+m^2)m^3n$.

    Write $C_k'=(D')^{-1}(\Id_m\otimes C_k^{\circ})D'=\bigoplus_{i=1}^m(D_i^{-1}C_k^{\circ}D_i)\in\M[n'](F)$. The above calculations show that
    \begin{align*}
        \Rank(\widehat{A_k'}-\widehat{C_k'}) & \le \Rank(\widehat{A_k'}-\widehat{C_k}) + \Rank(\widehat{C_k}-\widehat{C_k'}) \\
        & = \Rank(\widehat{A_k'}-\widehat{C_k}) + \Rank(C_k-C_k') \\
        & \le \Rank(\widehat{A_k'}-\widehat{C_k}) + \Rank(D^{-1}(\Id_m\otimes C_k^{\circ})D - (D')^{-1}(\Id_m\otimes C_k^{\circ})D') \\
        & \le \Rank(\widehat{A_k'}-\widehat{C_k}) + \Rank(D'-D) + \Rank((D')^{-1}-D^{-1}) \\
        & \le \Rank(\widehat{A_k'}-\widehat{C_k}) + 2\Rank(D'-D) \hspace{8em}\text{(by \Lref{lem:inv-diff})} \\
        & \le \eps'mn + 28\eps'(d+m^2)m^3n \le 29\eps'(d+m^2)m^3n.
    \end{align*}

    Now $D_1^{-1}C_1^{\circ}D_1,\dots,D_1^{-1}C_d^{\circ}D_1$ is an exact solution of $P_1,\dots,P_r$ since $C_1^{\circ},\dots,C_d^{\circ}$ is. Moreover, to show that it well-approximates $A_1,\dots,A_d$, let $n_0=\min\set{n,q}$ and $n_1=\max\set{n,q}$. The submatrix $X_0$ of $\widehat{A_k'}-\widehat{C_k'}$ containing its first $n_0$ rows and columns is the same as the submatrix of $\widehat{A_k}-\widehat{D_1^{-1}C_k^{\circ}D_1}$ containing its first $n_0$ rows and columns. Therefore
    \begin{align*}
        \Rank(\widehat{A_k}-\widehat{D_1^{-1}C_k^{\circ}D_1}) \le \Rank(X_0) + 2(n_1-n_0) &\le 29\eps'(d+m^2)m^3n+2\eps'(d+m^2)n \\
        &\le 31\eps'(d+m^2)m^3n=\eps n.
    \end{align*}
    This proves the rank-stability of $\Alg$.
\end{proof}

\section{Instability} \label{sec:instability}

In general, it is a hard task to prove that a given ``reasonable'' group is not stable. In \cite{BeckerLubotzkyThom19}, it was shown that Abels' group\mdash{}which is finitely presented, solvable and residually finite\mdash{}is not P-stable. In \cite{ElekGrabowski21}, it was shown that this group is also not rank-stable.
In this section we develop an instability criterion for algebras, inspired by the aforementioned results. This criterion is then applied to finite-dimensional Lie algebras, via their universal enveloping algebras.

We first recall the notion of amenable algebras (see \cite{Elek03,Gromov99}):
\begin{definition}
  Let $\Alg$ be a finitely generated algebra, and let $S$ be a finite set of generators. We say that $\Alg$ is \textbf{amenable} if for every $\eps>0$ there exists a finite-dimensional subspace $V\le \Alg$ such that
  \[
    \dim_F \left(SV+V\right) < (1+\eps)\cdot \dim_F V.
  \]
\end{definition}
Every algebra of subexponential growth is amenable \cite{Elek03}, and a group is amenable if and only if its group algebra is amenable \cite{Bartholdi08}. Notice that the above definition in fact refers to $\Alg$ as an amenable left module over itself; similarly, one can define right amenability.

\begin{theorem} \label{thm:unstable}
  Let $\Alg$ be a finitely presented algebra with a non-zero ideal $I\vartriangleleft \Alg$ such that:
  \begin{itemize}
      \item The ideal $I$ is finitely generated as a left ideal of $\Alg$;
      \item the algebra $\Alg/I$ is amenable; and
      \item $\Alg/I$ has no finite-dimensional representations.
  \end{itemize}
  Then $\Alg$ is not rank-stable.
\end{theorem}

\begin{remark}{\ }
\begin{enumerate}
    \item All of the required assumptions of Theorem \ref{thm:unstable} except for the last one are automatically satisfied for any ideal $I$ whenever $\Alg$ is a graded (left)  Noetherian algebra with finite-dimensional homogeneous components; such an algebra is always finitely presented \cite[Theorem~17]{Lewin74}, \cite[Theorem~2.2]{Greenfeld23} and residually finite-dimensional. Furthermore, every left ideal of a Noetherian ring is finitely generated, and by \cite{StephensonZhang97}, such an algebra has a subexponential growth, hence all of its homomorphic images are amenable.

    \item The assumption in Theorem \ref{thm:unstable} that $I \vartriangleleft \Alg$ is finitely generated as a one-sided ideal is crucial. Indeed, the free algebra $\Alg=F\left<x,y\right>$ admits a finitely generated ideal $I=\left<xy-yx-1\right>$ which satisfies all of the other assumptions; however, the free algebra is clearly rank-stable.
\end{enumerate}
\end{remark}

\begin{proof}[{{Proof of \Tref{thm:unstable}}}]
  Let $\Alg$ and $I$ be as in the statement of the theorem and assume to the contrary that $\Alg$ is rank-stable. Fix a presentation
  \[
    \Alg = F\left<x_1,\dots,x_d\right>/\left<P_1,\dots,P_r\right>,
  \]
  and let $S=F+Fx_1+\cdots+Fx_d$ be the standard generating subspace of $\Alg$. Let $m$ be the maximal degree of all of the $P_i$'s, and let $T$ be a basis of monomials for $S^m$. 
  By assumption, the ideal $I\neq 0$ is finitely generated as a left ideal, so we can write:
  \[
    I=\Alg f_1(x_1,\dots,x_d)+\cdots+\Alg f_\nu(x_1,\dots,x_d)
  \]

  for some $\nu\in\mathbb{N}$. Each $f_\mu$ (for $1\le \mu \le \nu$) is a polynomial in $x_1,\dots,x_d$ composed of at most, say, $l$ monomials; enlarging $m$ if necessary, we further assume that each one of these monomials is of length at most $m$.

  By the amenability of $\Alg/I$ we can take $(V_i)_{i=1}^{\infty}$ to be a sequence of finite-dimensional subspaces of $\Alg/I$ such that $TV_i\subseteq V_i + \mathcal{E}_i$, where
  \[
    \frac{\dim_F \mathcal{E}_i}{\dim_F V_i}\xrightarrow{i\rightarrow \infty} 0
  \]
  (here $T$ is identified with its image modulo $I$). Observe that $\dim_F V_i\xrightarrow{i\rightarrow \infty} \infty$, since by assumption $\Alg/I$ has no finite-dimensional representations and in particular there are no finite-dimensional subspaces that are invariant under $T$. Let
  \[
    V_i^{\circ} = \set{v\in V_i\ |\ \forall s\in T,\ sv\in V_i},
  \] be the ``$S^m$-interior'' of $V_i$, and let
  \[
    \psi\colon V_i/V_i^{\circ} \rightarrow \bigoplus_{s\in T} (V_i+\mathcal{E}_i)/V_i
  \]
  be the linear map defined by:
  \[
  \psi\colon  u+V_i^{\circ} \mapsto \left(su+V_i\right)_{s\in T}.
  \]
  By the definition of $V_i^{\circ}$, the map $\psi$ is well-defined and injective, hence
  \begin{equation}\label{eq:Vi-Vicirc}
    \dim_F V_i/V_i^{\circ} \le |T|\cdot \dim_F \mathcal{E}_i.
  \end{equation}
  We also let
  \[
    U_i = \{v\in V_i\ |\ \forall 1\le j\le d,\ x_jv\in V_i\}
  \] 
  be the ``$S$-interior'' of $V_i$.
  Clearly, $V_i^{\circ}\le U_i\le V_i$. Fix a vector space  complement $V_i=U_i\oplus W_i$.

Denote $n_i=\dim_F V_i$. For each $1\le j\le d$, we define a matrix $A_{i,j}\in\M[n_i](F)$ representing the linear endomorphism of $V_i$ whose action on $U_i$ is given by $u\mapsto x_ju$ and zero on $W_i$. We claim that $\rank(P_j(A_{i,1},\dots,A_{i,d}))\xrightarrow{i\to\infty} 0$ for all $j$, while for sufficiently large $i$, these matrices are not well-approximated by any exact solution to $P_1,\dots,P_r$.

  Indeed, fix a relation $P\in\set{P_1,\dots,P_r}$ and write
  \[
    P(\vec{x})=\sum_{\underline{\alpha}=(\alpha_1,\dots,\alpha_p)} c_{\underline{\alpha}} x_{\alpha_1}\cdots x_{\alpha_p}
  \]
  (each $p$ depends on $\underline{\alpha}$, but recall that $p\le m$).
  Now for each $u\in V_i^{\circ}$ and for each $\underline{\alpha}$ observe that:
  \begin{align*}
  x_{\alpha_p}u,x_{\alpha_{p-1}}x_{\alpha_p}u,\dots,x_{\alpha_1}\cdots x_{\alpha_p}u &\in V_i \\
  u,x_{\alpha_p}u,\dots,x_{\alpha_2}\cdots x_{\alpha_p}u &\in U_i,
  \end{align*}
  so
  \begin{align}\label{eqs:PAi}
    P(A_{i,1},\dots,A_{i,d})u & = \sum_{\underline{\alpha}} c_{\underline{\alpha}} A_{i,\alpha_1}\cdots A_{i,\alpha_p} u\\\notag
    & = \sum_{\underline{\alpha}} c_{\underline{\alpha}} A_{\alpha_1}\cdots A_{\alpha_{p-1}} x_{\alpha_p} u\\\notag
    & \quad\vdots & \\ \notag
    & = \sum_{\underline{\alpha}} c_{\underline{\alpha}} x_{\alpha_1}\cdots x_{\alpha_p} u = 0.
  \end{align}
  Thus, $V_i^{\circ}\subseteq \ker P(A_{i,1},\dots,A_{i,d})$ and consequently, by $\eqref{eq:Vi-Vicirc}$,
  \[
    \rank(P(A_{i,1},\dots,A_{i,d}))\le \frac{1}{n_i}\dim_F V_i/V_i^{\circ} \le \frac{1}{n_i}|T|\dim_F \mathcal{E}_i =: \omega_i \xrightarrow{i\rightarrow \infty} 0.
  \]

  Recall that we are assuming that $\Alg$ is rank-stable.
  Let $\eps=\min\set{\frac{1}{4\nu l m},\frac{1}{4d\nu}}$, and choose $\delta>0$ given from the rank-stability of $\Alg$ for this choice of $\eps$. Now choose $i\gg 1$ such that $n_i>\frac{4}{1-4\eps d}$ and $\omega_i < \min\set{\frac{1}{4\nu},\delta}$. (From now on we omit the $i$ subscript of $A_{i,j}$.) By rank-stability, there exists $n_i'$ and a tuple of $n_i'\times n_i'$ matrices $B_1,\dots,B_d\in\M[n_i'](F)$ such that $P_j(\vec{B})=0$ for each $1\le j\le r$ and $\Rank(\widehat{A_e}-\widehat{B_e})<\eps n_i$ for each $1\le e \le d$. Also, by \Lref{lem:small-approx}, by taking $\delta<(\eps d)^2$ we may assume that
  \[
    (1-\eps d)n_i\le n_i'\le (1+\eps d)n_i.
  \]

  Recall that $I=\Alg f_1(x_1,\dots,x_d)+\cdots+\Alg f_\nu(x_1,\dots,x_d)$.
  By the same argument as \eqref{eqs:PAi}, we get that for each $1\le \mu\le \nu$ (recalling that $f_\mu$ is a sum of monomials of degrees at most $m$), if $u\in V_i^{\circ}$, then 
  \[
    f_\mu(A_1,\dots,A_d)u=f_\mu(x_1,\dots,x_d)u=0.
  \]
  So by $\eqref{eq:Vi-Vicirc}$, we have that
  \[
  \Rank(f_\mu(A_1,\dots,A_d))\le \dim_F V_i/V_i^{\circ} \leq n_i \omega_i<\frac{n_i}{4\nu}
  \]
  for each $1\le \mu \le \nu$.
  Hence by \Lref{lem:computing polyrank},
  \begin{align*}
  \Rank(f_\mu(B_1,\dots,B_d)) & <  \frac{n_i}{4\nu} + l m\eps n_i + \left|n_i-n_i'\right| \\
  & \le \frac{n_i}{4\nu} + \frac{n_i}{4\nu} + \eps dn_i \\
  &  \le \frac{3n_i}{4\nu}.    
  \end{align*}
  Therefore, the subspace
  \[
    Z \coloneqq \bigcap_{\mu=1}^{\nu} \ker f_\mu(B_1,\dots,B_d) \leq F^{n'_i}
  \]
  satisfies
  \[
    \dim_F Z \ge n'_i-\nu \cdot \frac{3}{4\nu}\cdot n_i \geq (1-\eps d)n_i -  \nu \cdot \frac{3}{4\nu}\cdot n_i = \frac{1-4 \eps d}{4} \cdot n_i \geq 1.
  \]
  We thus get a non-trivial subspace $Z\le F^{n'_i}$ on which $\Alg$ acts via $B_1,\dots,B_d$. Indeed, fix $1\le \mu_0\le\nu$, and let $1\le \rho\le d$. Since $I\vartriangleleft \Alg$, in the algebra $\Alg$ we can write
  \[
    f_{\mu_0}(x_1,\dots,x_d)x_{\rho}=\sum_{\mu=1}^{\nu}w_{\mu}f_{\mu}(x_1,\dots,x_d)
  \]
  where $w_\mu$ are some non-commutative polynomials in $x_1,\dots,x_d$. Since $P_j(\vec{B})=0$ for all $j$, we may substitute $x_1=B_1,\dots,x_d=B_d$, and conclude that for each $1\leq \rho\leq d$ and $v\in Z$:
  \[
    f_{\mu_0}(B_1,\dots,B_d)\cdot B_\rho v = \sum_{\mu=1}^\nu w_\mu f_\mu(B_1,\dots,B_d) \cdot v = 0
  \]
  so $B_\rho v\in Z$, proving that $\Alg$ acts on $Z$ via $x_i\cdot v=B_i\cdot v$ for each $1\leq i\leq d$.
  Moreover, notice that $I$ annihilates $Z$ under this action, so we obtain a representation of $\Alg/I$ on $Z$, contradicting the assumption that $\Alg/I$ has no finite-dimensional representations.
\end{proof}

\Tref{thm:unstable} applies for the projectivization of any finitely presented amenable simple algebra. For instance:

\begin{corollary} \label{cor:Weyldef}
  Let
  \[
    \overline{\Alg}_1=\C\left<x,y,t\right>/\left<xt-tx,\ yt-ty,\ xy-yx-t^2 \right>
  \]
  be the projectivization of the first Weyl algebra $\Alg_1$. Then $\overline{\Alg}_1$ is a finitely presented, residually finite-dimensional non rank-stable algebra.
\end{corollary}

\begin{proof}
Notice that $\overline{\Alg}_1$ is residually finite-dimensional, being graded with finite-dimensional homogeneous components. The ideal $\left<t-1\right>\vartriangleleft \overline{\Alg}_1$ is a principal left ideal, and the quotient $\overline{\Alg}_1/{\left<t-1\right>}\cong \Alg_1$ is infinite-dimensional and simple, hence has no finite-dimensional representations, and also amenable, since it has polynomial growth. By \Tref{thm:unstable}, the algebra $\overline{\Alg}_1$ is non rank-stable.
\end{proof}

Unlike associative algebras, finite-dimensional Lie algebras tend to be non rank-stable:

\begin{theorem} \label{thm:Lie}
Let $\mathfrak{g}$ be a finite-dimensional Lie algebra over $\C$. If $\mathfrak{g}$ is non-solvable, or admits a non-abelian nilpotent homomorphic image, then it is non rank-stable.
\end{theorem}

\begin{proof}
The universal enveloping algebra $U(\mathfrak{g})$ is finitely presented and Noetherian; furthermore, it has polynomial growth by the Poincar\'e--Birkhoff--Witt theorem (in fact, its Gel'fand--Kirillov dimension is $\GK U(\mathfrak{g}) = \dim_\C \mathfrak{g}$). Therefore, to apply \Tref{thm:unstable}, it suffices to show that $U(\mathfrak{g})$ admits an ideal not contained in any finite-codimensional ideal, or even an infinite-codimensional maximal ideal.

If $\mathfrak{g}$ admits a non-abelian nilpotent homomorphic image $\mathfrak{n}$ then a theorem of Dixmier~\cite[Theorem~4.7.9]{Dixmier96} ensures $U(\mathfrak{g})$ admits an $n$-th Weyl algebra $\mathcal{A}_n\cong \mathcal{A}_1^{\otimes n}$, which is simple and infinite-dimensional, as a homomorphic image. 
Here is a more elementary and straightforward argument. 
If $\mathfrak{g}$ admits a non-abelian nilpotent homomorphic image $\mathfrak{n}$ then we may assume that $[[\mathfrak{n},\mathfrak{n}],\mathfrak{n}]=0$ and pick $x,y\in \mathfrak{n}$ such that $0\neq z=[x,y]\in [\mathfrak{n},\mathfrak{n}]\subseteq Z(U(\mathfrak{n}))$. Consider the quotient $U(\mathfrak{n})/(z-\alpha)$ for some $\alpha$, in which $\alpha^{-1} x,y$ generate a copy of the first Weyl algebra $\mathcal{A}_1$. We found a homomorphic image of $U(\mathfrak{g})$ with no finite-dimensional homomorphic images.

If $\mathfrak{g}$ is not solvable, we may replace it with a simple homomorphic image of itself.
Let $\chi\colon Z(U(\mathfrak{g}))\rightarrow \mathbb{C}$ be a generic central character and let $I$ be the ideal of $U(\mathfrak{g})$ generated by $\ker(\chi)$.
Then $I$ is equal to the annihilator $\ann(M)$ of a Verma module $M$, by \cite{Joseph98}. Since $\chi$ is generic, the module $M$ is infinite-dimensional and hence $I$ is infinite-codimensional.

We claim that $I$ is a maximal ideal of $U(\mathfrak{g})$. Indeed, pick $I\subseteq J$ a maximal (in particular, primitive) ideal containing it. By \cite{Duflo77}, $J=\ann(L)$ for some highest weight module $L$. Again, since $\chi$ is generic, all Verma modules with this central character are irreducible, so $J$ is in fact the annihilator of a Verma module. Hence $J$ is centrally generated (see \cite{Joseph98}), namely, $J$ is generated by $J\cap Z(U(\mathfrak{g}))=I\cap Z(U(\mathfrak{g}))$, so $I=J$ and $I$ is maximal. The proof is completed.
\end{proof}

For an abelian Lie algebra $\mathfrak{a}$ of dimension greater than one, the universal enveloping algebra $U(\mathfrak{a})$ is a polynomial ring in several variables, so its rank-stability is equivalent to the rank-stability of a polynomial ring.

\smallskip

\noindent \textit{Acknowledgements.} 
We thank Maria Gorelik, Boris Kunyavskii and Shifra Reif for helpful discussions. The second author is supported by the Bar-Ilan President's Doctoral Fellowships of Excellence and by Israeli Science Foundation grant \#957/20.

\bibliographystyle{abbrv}
\bibliography{rank}

\end{document}